\documentclass[12pt]{amsart}
\usepackage{fullpage}
\usepackage{amsfonts}
\usepackage{amssymb}
\usepackage{amscd}
\usepackage{latexsym}

\def\SBIMSMark#1#2#3{
 \font\SBF=cmss10 at 10 true pt
 \font\SBI=cmssi10 at 10 true pt
 \setbox0=\hbox{\SBF Stony Brook IMS Preprint \##1}
 \setbox2=\hbox to \wd0{\hfil \SBI #2}
 \setbox4=\hbox to \wd0{\hfil \SBI #3}
 \setbox6=\hbox to \wd0{\hss
             \vbox{\hsize=\wd0 \parskip=0pt \baselineskip=10 true pt
                   \copy0 \break%
                   \copy2 \break%
                   \copy4 \break}}
 \dimen0=\ht6   \advance\dimen0 by \vsize \advance\dimen0 by 8 true pt
                \advance\dimen0 by -\pagetotal
 \dimen2=\hsize \advance\dimen2 by .25 true in
  \openin2=publishd.tex
  \ifeof2\setbox0=\hbox to 0pt{}
  \else 
     \setbox0=\hbox to 3.1 true in{
                \vbox to \ht6{\hsize=3 true in \parskip=0pt  \noindent  
                \input publishd.tex 
                \vfill}}
  \fi
  \closein2
  \ht0=0pt \dp0=0pt
 \ht6=0pt \dp6=0pt
 \setbox8=\vbox to \dimen0{\vfill \hbox to \dimen2{\copy0 \hss \copy6}}
 \ht8=0pt \dp8=0pt \wd8=0pt
 \copy8
 \message{*** Stony Brook IMS Preprint #1, #2 ***}
}

\catcode`\@=11

\long\def\@savemarbox#1#2{\global\setbox#1\vtop{\hsize\marginparwidth 
  \@parboxrestore\tiny\raggedright #2}}
\catcode`\@=12

\def\Empty{}
\newcommand\oplabel[1]{
  \def\OpArg{#1} \ifx \OpArg\Empty {} \else
  	\label{#1}
  \fi}
	
\long\def\realfig#1#2#3#4{
\begin{figure}[htbp]
\centerline{\psfig{figure=#2,width=#4}}
\caption[#1]{#3}
\oplabel{#1}
\end{figure}}

\newcommand{\comm}[1]{}

\usepackage{psfig}

\newtheorem{thm}{Theorem}[section]
\newtheorem{conjecture}{Conjecture}[section]

\newtheorem{lem}[thm]{Lemma}
\newtheorem{prop}[thm]{Proposition}
 
\newtheorem{yoc}{Yoccoz Inequality}
\newtheorem{rem}[thm]{Remark}

\newenvironment{pf}{\proof[\proofname]}{\endproof}
\newenvironment{pf*}[1]{\proof[#1]}{\endproof}

\theoremstyle{definition}

\usepackage[OT2,OT1]{fontenc}
\usepackage{euscript}
\def\cyr{\fontencoding{OT2}\fontfamily{wncyr}\selectfont}

\newcommand{\cal}[1]{{\mathcal #1}}
\newcommand{\BBB}[1]{{\mathbb #1}}

\newcommand{\n}{\operatorname{and}}

\newcommand{\dist}{\operatorname{dist}}

\newcommand{\Per}{\operatorname{Per}}
\newcommand{\h}{{\cal H}}
\newcommand{\M}{{\cal M}}
\newcommand{\rt}{{\operatorname{root}}}

\newcommand{\he}{\underset{\operatorname{hb}}{\sim}}

\newcommand{\ap}{{\cal A}}
\newcommand{\rp}{{\cal R}}

\renewcommand{\mod}{\operatorname{mod}}

\newcommand{\tl}{\tilde}

\newcommand{\eps}{\epsilon}

\numberwithin{equation}{section}

\newcommand{\thmref}[1]{Theorem~\ref{#1}}
\newcommand{\propref}[1]{Proposition~\ref{#1}}
\newcommand{\secref}[1]{\S\ref{#1}}
\newcommand{\lemref}[1]{Lemma~\ref{#1}}
 
\newcommand{\figref}[1]{Fig.~\ref{#1}}

\newcommand{\intr}{\curlyvee}

\newcommand{\Cub}{{\EuScript C}}

\newcommand{\CubR}{{\EuScript C}_{\BBB R}}

\newcommand{\CubM}{{\EuScript C}^{\#}}

\newcommand{\CubMR}{{\EuScript C}^{\#}_{\BBB R}}

\newcommand{\sector}[2]{\backslash #1,#2/}

\newcommand{\lren}{{\operatorname{R}}}

\newcommand{\rren}{\cyr{Ya}}

\newcommand{\lrren}{\lren\times\rren}

\renewcommand{\Re}{\operatorname{Re}}

\newcommand{\const}{\operatorname{const}}
\begin{document}
\bibliographystyle{amsalpha}

\title[Intertwining Surgery]{Geography of the Cubic Connectedness Locus I:
Intertwining Surgery}

\author{Adam Epstein}
\address{Institute for Mathematical Sciences\\
Department of Mathematics\\
State University of New York at Stony Brook\\
New York, 11794-3660}
\email{adame@@math.sunysb.edu}

\author{Michael Yampolsky}
\address{Department of Mathematics\\
State University of New York at Stony Brook\\
New York, 11794-3651}
\email{yampol@@math.sunysb.edu}

\maketitle
\SBIMSMark{1996/10}{August 1996}{}
\tableofcontents
\newpage
\section{Introduction}

The prevalence of Mandelbrot sets in one-parameter complex analytic 
families is a well-studied phenomenon in conformal dynamics.
Its explanation gave rise to the theory of {\it renormalization} \cite{DH},
and subsequent efforts to invert this procedure by means of
{\it surgery} on quadratic polynomials \cite{BD,BF}. 

In this paper we  exhibit  products
of Mandelbrot sets in the two-dimensional complex parameter space
of cubic polynomials. These products were observed by J. Milnor
in computer experiments which inspired Lavaurs'
proof of non local-connectivity for the cubic connectedness locus
\cite{Lav}. Cubic polynomials 
in such a product may be renormalized to produce a pair of
quadratic maps.  The inverse construction is an  {\it intertwining surgery} on
 two quadratics. The idea of intertwining first appeared in
a collection of problems edited by Bielefeld \cite{Biel}.
Using quasiconformal surgery techniques of Branner and Douady \cite{BD},
we show that any two quadratics may be intertwined to obtain
a cubic polynomial. The proof of continuity in our two-parameter setting
requires further considerations involving ray combinatorics and a 
pullback argument.

After this project was finished, we were informed by P. Haissinsky
that he is independently working on related problems \cite{Hai}.

\medskip
\noindent
{\bf Acknowledgments:}
This paper was motivated by J. Milnor's Autumn 1995 Stony Brook
 lectures on the dynamics
of cubic polynomials, and
developed out of joint meditation of the two authors in front of the
the full-color version of \figref{AB-plane}. We thank J. Milnor
for numerous discussions of our results
and many helpful suggestions as this paper progressed.
We are indebted to M. Lyubich for fruitful conversations
concerning various aspects of quadratic dynamics. 
Further thanks are due to J. Kiwi for sharing his understanding of cubic
maps and discussing some of his current work, and to X.~Buff for
communicating his results. This project was conducted in
the congenial  atmosphere of IMS at Stony Brook, and we thank
our colleagues for their interest and moral support.

The computer pictures in this paper were produced using software
written by J.~Milnor and S. Sutherland.

\section{Preliminaries }
\label{preliminaries}

In this section we discuss the relevant facts and tools of holomorphic
dynamics. We assume that the reader is familiar with the basic notions and
principles of the theory of quasiconformal maps 
(see \cite{Le} for a comprehensive
account). The  knowledgeable reader is invited to proceed directly
to \secref{outline}.

\subsection{Polynomial dynamics}
{\bf Julia sets, external rays, landing theorems,
combinatorial rotation number, Yoccoz Inequality}

We recall the basic definitions and results in the theory
of polynomial dynamics. Supporting details may be found in
\cite{Mil-survey}.

Let $P:\BBB C\to \BBB C$ be a complex polynomial of degree $d\ge 2$.
The {\it filled Julia set} of $P$
is defined as
$$ K(P)=\{z\in \BBB C|\{P^{\circ n}(z)\} \operatorname{is\; bounded}\}$$
and the {\it Julia set} as $J(P)=\partial K(P)$.
Both of these are nonempty compact sets which are connected
if and only all critical points of $P$ have bounded orbits.

Recall that if $P$ is a monic polynomial with connected Julia set
then there exists a unique analytic homeomorphism (the 
{\it B\"ottcher map})
$$B_P:\BBB C\setminus K(P)\to \BBB C\setminus \bar{\BBB D}$$
which is tangent to the identity at infinity, that is  $B_P(z)/z\to 1$
 as $z\to \infty$.
The B\"ottcher map conjugates $P$ to $z\mapsto z^d$,
$$B_P(P(z))=(B_P(z))^d,$$
thereby determining a dynamically natural polar coordinate system on
$\BBB C\setminus K(P)$.
For $\rho>1$ the {\it equipotential} $E_\rho$ is the inverse image under $B_P$
of the circle $\{\rho e^{2\pi i\theta}|\theta\in{\BBB R}\}$.
The {\it external ray at angle $\theta$} is similarly defined as the 
inverse image $r_\theta$ of the radial line 
$\{ \rho e^{2\pi i\theta}|\rho>1\}$.
Since $P$ maps $r_\theta$ to $r_{d\cdot \theta}$,
the  ray $r_\theta$ is {\it periodic} if and only if the angle $\theta$ is
periodic (mod $1$) under multiplication by $d$.
An  external ray $r_\theta$ is said to {\it land}
at a point $\zeta\in J(P)$ when 
$$ \lim_{\rho\to 1}B_P^{-1}(\rho e^{2\pi i\theta})=\zeta.$$
We note that if  the Julia set of $P$ is locally connected then all 
rays $r_\theta$ land, and their endpoints  depend continuously on the
angle $\theta$ (see the discussion in \cite{Mil-survey}). We refer to 
\cite{Mil-survey} for the proofs of  the following results:

\begin{thm}[Sullivan, Douady and Hubbard] 
\label{landing-1}
If K(P) is connected, then every
periodic external ray lands at a periodic point which is either repelling
or parabolic.
\end{thm}

\begin{thm}[Douady, Yoccoz]
\label{landing-2}
If $K(P)$ is connected, every repelling or parabolic periodic point is 
the landing point of at least one external ray which is necessarily periodic.
\end{thm}

The landing points of such rays depend continuously on parameters:

\begin{lem}[\cite{GM}]
\label{same-rays}
Let $P_t$ be a continuous family of monic degree $d$ polynomials
with continuously chosen repelling periodic points $\zeta_t$. If the
ray of angle $\theta$ for $P_{t_0}$ lands at $\zeta_{t_0}$, then for all
$t$ close to $t_0$ the ray of angle $\theta$ for $P_t$ lands at $\zeta_t$.
\end{lem}

Kiwi has proved the following useful separation principle which directly
illustrates why a degree $d$ polynomial admits at most $d-1$ non-repelling
periodic orbits; the latter result was earlier shown by Douady and Hubbard
and appropriately generalized to rational maps by Shishikura.

\begin{thm}
\label{separation} 
Let $P$ be a polynomial with connected Julia set, $n$ a common multiple
of the periods of non-repelling periodic points, $\cal R$ the union of all
external rays fixed under $P^{\circ n}$ together with their landing points,
and $U_1,\ldots,U_m$ be the connected components of 
$\BBB C\setminus\bigcup_{j=0}^n P^{-\circ j}(\cal R)$. Then:
\begin{itemize}
\item Each component $U_i$ contains at most one non-repelling periodic point;
\item Given any non-repelling periodic orbit $\zeta_1,\ldots,\zeta_\ell$ 
passing through $U_{i_1},\ldots,U_{i_\ell}$, at least one of the
components $U_{i_k}$ also contains some critical point.
\end{itemize}
\end{thm}

We assume henceforth that $K(P)$ is connected. 
Let $r=r_\theta$ be a periodic external ray landing at the periodic
point $\zeta\in K(P)$, whose orbit we enumerate
$$\zeta=\zeta_0\mapsto \zeta_1\mapsto\ldots\mapsto\zeta_n=\zeta.$$
Denote by $A_i\subset \BBB Q/\BBB Z$ the set of angles of the rays in
the orbit of $r$ landing at $\zeta_i$.
The iterate $P^{\circ n}$ fixes each point $\zeta_i$, permuting 
the various rays landing there while  preserving their cyclic order.
Equivalently, multiplication by $d^n$ carries the set $A_i$
onto itself by an order-preserving bijection. For each
$i$ we may label the angles in $A_i$ as 
$0\leq \theta^1<\theta^2<\ldots<\theta^q<1$; then
$$d^n\theta^i\equiv\theta^{i+p}\;\;(\mbox{mod } 1)$$
for some integer $p$, and we refer to the ratio $p/q$ as the
{\it combinatorial rotation number} of $r$. 
The following theorem of Yoccoz (see \cite{Hub}) relates the 
combinatorial rotation number of a ray landing at a period $n$ point
$\zeta$ to the multiplier $\lambda=(P^{\circ n})'(\zeta)$.
\begin{yoc}
Let $P$ be a monic polynomial
with  connected Julia set, and  $\zeta\in K(P)$  a repelling fixed
point  with multiplier $\lambda$.
If $\zeta$ is  the landing point of $m$ distinct cycles of external rays
with  combinatorial rotation number $p/q$ then 
\begin{equation}
\label{yoccoz}
\frac{\operatorname{Re}\rho}{|\rho-2\pi i p/q|^2}\geq 
  \frac{mq}{2\log d}.
\end{equation}
where $\rho$ is the suitable choice of $\log\lambda$.
\end{yoc}
More geometrically, the inequality asserts that $\rho$  lies in 
the closed disc of radius $\log d/(mq)$ tangent to the imaginary axis
at $2\pi i p/q$. 

\subsection{Polynomial-like maps}
{\bf Hybrid equivalence, Straightening Theorem,
continuity of straightening}

Polynomial-like mappings, introduced by Douady and Hubbard in \cite{DH},
are a key tool in holomorphic dynamics. A
{\it polynomial-like mapping}  {\it of degree $d$} is
a proper degree $d$ holomorphic map $f:U\to V$  between topological discs,
where $U$ is compactly contained in $V$.
One  defines the filled Julia set 
$$K(f)=\{z\in U| f^{\circ n}(z)\in U,\;\forall n\geq 1\}$$
and the Julia set
$J(f)=\partial K(f)$.
We say that the map $f$ is {\it quadratic-like} if the degree $d=2$,
and {\it cubic-like} if $d=3$.

Polynomial-like maps $f:U\to V$ and $\tl f:\tl U\to \tl V$
are {\it hybrid equivalent}
$$f\he \tl f$$
if there exists a quasiconformal
homeomorphism $h$ from a neighborhood of $K(f)$
to a neighborhood of $K(\tl f)$, such that $h\circ f=\tl f\circ h$
near $K(f)$ and $\bar\partial h=0$ almost everywhere on $K(f)$.
We remark that   $h$  can be chosen to
be a conjugacy between $f|_U$ and $\tl f|_{\tl U}$. Notice that
$h$ is conformal on the interior of $K(f)$ and therefore preserves the 
multipliers of attracting periodic orbits. In view of the well-known 
quasiconformal invariance of indifferent multipliers, we observe:

\begin{rem} \label{eigrmk}
A hybrid equivalence between polynomial-like maps sends repelling
to repelling orbits, and preserves the multipliers of attracting and
indifferent orbits.
\end{rem}

The following is fundamental:

\begin{thm}[Straightening Theorem, \cite{DH}]
\label{strt-thm}
Every polynomial-like mapping $f:U\to V$ of degree $d$ is hybrid equivalent
to a polynomial $P$ of degree $d$. If $K(f)$ is connected then $P$ is unique
up to conjugation by an affine map. 
\end{thm}

\noindent
For quadratic-like $f$ with connected Julia set, we write
$\chi(f)=c$ where $$f_c(z)=z^2+c$$ is the unique 
hybrid equivalent polynomial.
The following theorem is due to Douady
and Hubbard; we employ the formulation of \cite[Prop. 4.7]{McM2}:
\begin{thm}
\label{straightening}
Let $f_k:U_k\to V_k$ be a sequence of quadratic-like maps with
connected Julia sets, which 
converges uniformly to a quadratic-like map $f:U\to V$ on a neighborhood
of $K(f)$. Then $\chi(f_k)\to\chi(f).$
\end{thm}

The proof of the uniqueness assertion in \thmref{strt-thm}
relies essentially on the following general lemma due to Bers \cite{Le}:

\begin{lem}
\label{rickman}
Let $U\subset \BBB C$ be open, $K\subset U$ be compact, and $\phi$ and $\Phi$
be two mappings $U\to\BBB C$ which are homeomorphisms onto their images.
Suppose that $\phi$ is quasiconformal, that $\Phi$ is quasiconformal
on $U\setminus K$, and that $\phi=\Phi$ on $K$. Then $\Phi$ is
quasiconformal, and $\bar\partial \phi=\bar \partial \Phi$ almost
everywhere on $K$.
\end{lem}

\subsection{Quadratic polynomials}
{\bf Mandelbrot set, renormalizable 
maps and tuning}
\label{quadratics}

Basic facts on the structure of the Mandelbrot set are found
in \cite{orsay-notes}. Our account of renormalization and the
Yoccoz construction follows 
\cite{Ly2}, see also \cite{Mil-rays}, and \cite{McM}.

The {\it connectedness locus} of the quadratic family $f_c(z)=z^2+c$ 
is the ever-popular Mandelbrot set 
$$\M=\{c\in \BBB C|\; J(f_c) \operatorname{is\; connected}\}$$
depicted in \figref{mandelbrot}. The following results  
are shown in \cite{orsay-notes}.

\begin{thm}[Douady and Hubbard]
 The Mandelbrot set is  compact and connected, with connected complement.
\end{thm}

By definition, the {\it hyperbolic components} of $\M$ are the
connected components $H$ of 
$\overset{\circ}{\M}$ such that $f_c$ has an attracting periodic orbit for 
$c\in H$. Recalling that there can be at most one such orbit, we denote
its multiplier $\lambda_H(c)$.

\begin{thm}[Douady and Hubbard]
Let $H$ be a hyperbolic component.
The multiplier map 
$$\lambda_H:H\to \BBB D$$
is a conformal isomorphism.
This map extends to a homeomorphism between $\bar H$ and the closed disc
$\bar{\BBB D}$.
\end{thm}

\realfig{mandelbrot}{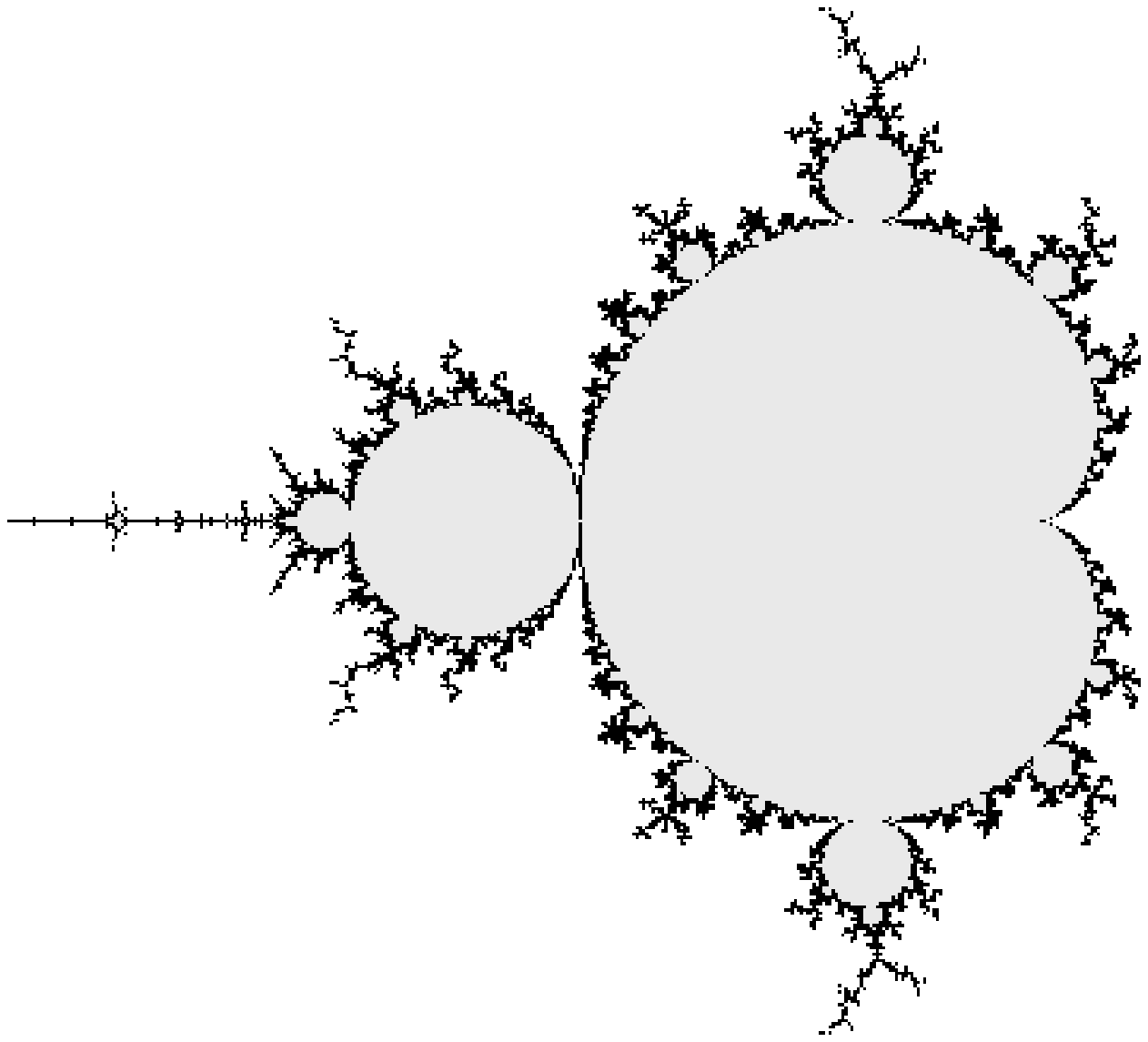}{The Mandelbrot set $\M$}{8cm}

Let $f_c$ be a quadratic polynomial with connected Julia set.
By  \thmref{landing-1} the external
ray of external argument $0$ lands at a fixed point of $f_c$, 
necessarily repelling or parabolic with multiplier 1, 
henceforth denoted $\beta_{f_c}$. The main
hyperbolic component $H_0$ is the set of all $c$ for which the other fixed
point $\alpha_{f_c}$ is attracting; the boundary point $c=1/4$ is hereafter
referred to as the {\it root} of $\M$.
 
The {\it $p/q$-limb} $L_{p/q}$ is the connected component of
$\M\setminus H_0$ whose boundary contains 
$$\rt_{p/q}=\lambda_{H_0}^{-1}(e^{2\pi i p/q})$$
and we denote by $H_{p/q}$ the hyperbolic component attached to $H_0$
at this point; it is always assumed that $(p,q)=1$. Notice that $L_{0/1}$ is
$\M$ itself. In view of the following, we may refer to $\alpha_{f_c}$ as
the {\it dividing fixed point}.

\begin{lem}
\label{rot-number}
For $q\geq 2$, a parameter value $c\in\M$ lies in $L_{p/q}$ if and only if
$\alpha_{f_c}$ is the landing point of an external ray
with combinatorial rotation number $p/q$.
\end{lem}

Consider a polynomial $f_c$ with connected Julia set.
Let $\zeta_1\mapsto \zeta_2\mapsto\ldots\mapsto\zeta_n=\zeta_1$ be a
repelling cycle of $f_c$, such that each $\zeta_i$ is the
landing point of at least two external rays. Let $\cal R$ be the 
collection of all external rays landing at these points,
and let $\cal R'=-\cal R$ be the symmetric collection.
Let us also choose an arbitrary equipotential $E$.
Denote by  $\Omega$ the component of 
$\BBB C\setminus (\cal R\cup \cal R'\cup E)$ containing $0$.
This region is bounded by four pieces of external rays and
two pieces of $E$.
Let $n$ be the period of these rays, $\zeta=\zeta_i$
the element of the cycle contained in $\partial \Omega$,
and $\Omega'\subset \Omega$ the 
component of $f_c^{\circ -n}(\Omega)$ attached to  $\zeta$.
If $0\in\Omega'$ then $f_c^{\circ n}:\Omega'\to \Omega$ is a branched
cover of degree 2.

Following Douady and Hubbard, we say that a  polynomial $f_c$ is 
{\it renormalizable} if there exists a repelling cycle $\{\zeta_i\}$
as above, such that $0\in\Omega'$ and 
$0$ does not escape $\Omega'$ under iteration of $f_c^{\circ n}$.
In this case $f_c^{\circ n}|_{\Omega'}$ can be extended to a 
quadratic-like map $f_c^{\circ n}:U\to V$ with connected Julia set
by a {\it thickening procedure} (a version of this procedure is employed in 
\secref{renormalization}). To emphasize the dependence of this construction
on the choice of periodic orbit, we shall say that 
this {\it renormalization} of $f_c$ is associated to  $\zeta$.

Recall that the $\omega$-{\it limit set} of a point $z$ under a map $f$
is defined as 
$$\omega_f(z)=\{w| f^{\circ n_k}(z)\to w \mbox{ for some } n_k\to\infty\}.
$$ When $f=f_c$ we simply write $\omega_c(z)$ and pay special attention to
the $\omega$-limit set of the critical point 0.
The following observation will be useful along the way:

\begin{rem}
\label{not-accumulate}
For a renormalizable quadratic polynomial $f_c$ with $n$ as above,
$$\omega_c(0)\subset \displaystyle \bigcup_{i=0}^{n-1} f_c^{\circ i}(\overline{\Omega'})\cap J(f_c).$$
In particular, $\beta_{f_c}\not\in\omega_c(0)$.
\end{rem}

\smallskip
\begin{thm}[Douady and Hubbard, \cite{DH}] 
Let  $f_{c_0}$ be a  renormalizable quad\-ratic polynomial with associated
periodic point $\zeta$.
Then  there exists a canonical embedding of the Mandelbrot set
$\M$ onto a subset $ \M'\ni c_0$ such that every map
$f_c$ with $c\in \M'\setminus\{\operatorname{one\:  point}\}$ 
is renormalizable with associated repelling periodic point $\zeta_c$,
where $c\mapsto \zeta_c$ is continuous and $\zeta_{c_0}=\zeta$.
\end{thm}

These subsets $\M'$ are customarily referred to as the {\it small copies}
of the Mandelbrot set. 
The inverse homeomorphism $\kappa:\M'\to\M$ is defined in terms of the
straightening map $\chi$:
$$\M'\ni c\;\;\mapsto\;\;  f_c^{\circ n}:U_c\to V_c
\;\;\overset{\chi}{\mapsto}\;\; \kappa(c)\in \M.$$
The periodic point $\zeta_c$ becomes parabolic with multiplier 1 
at the excluded parameter value, hereafter referred to as 
the {\it root} of $\M'$. 
We write $\M_{p/q}$ for the small copy ``growing'' from the hyperbolic
component $H_{p/q}$, its root being the point root$_{p/q}$.

\subsection{Cubic polynomials}{\bf The connectedness locus, types of
hyperbolic  components, Per$_n(\lambda)-$curves, real cubic family}

We now turn our attention to cubic polynomials. Our presentation follows the 
detailed discussion in \cite{Mil}. 

Observe that every cubic polynomial is
affine conjugate to a map of the form
\begin{equation}
\label{normal-form}
F_{a,b}(z)=z^3-3a^2z+b,
\end{equation}
with critical points $a$ and $-a$.
This normal form is unique up to conjugation by $z\mapsto -z$,
which interchanges $F_{a,b}$ and $F_{a,-b}$.
The pair of complex numbers $A=a^2$ and $B=b^2$ parameterize the space
of cubic polynomials modulo affine conjugacy.

The {\it cubic connectedness locus} is the set
$\Cub\subset \BBB C^2$ of all pairs $(A,B)$ for which the corresponding
polynomial $F_{a,b}$ has connected Julia set. As in the quadratic
case, the connectedness locus is compact and connected with connected
complement. These results were obtained by Branner and Hubbard \cite{BrHu}
who showed moreover that this set is {\em cellular}, the
intersection of a sequence of strictly nested closed discs. On the other hand,
Lavaurs \cite{Lav} proved that $\Cub$ is not locally
connected (compare Appendix \ref{section-comb}).

\realfig{AB-plane}{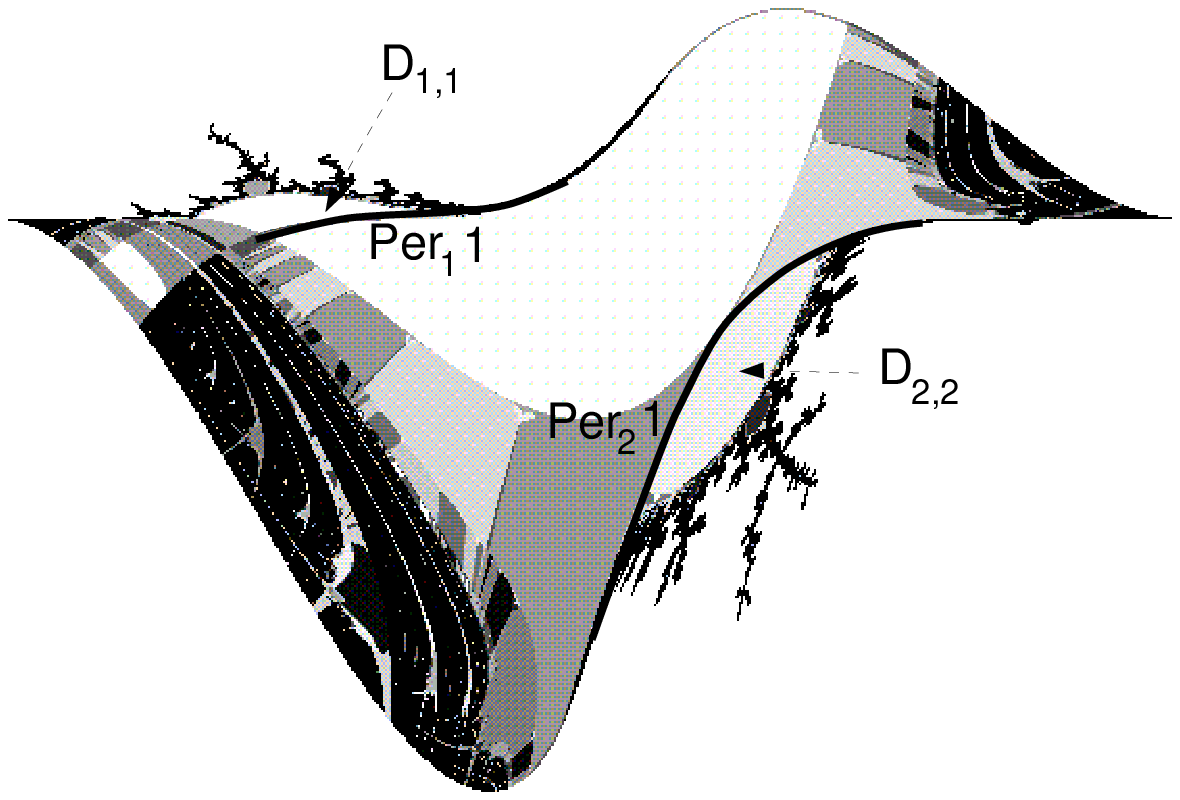}{Connectedness locus $\CubR$ in the real
cubic family}{11cm}

Milnor distinguishes four different types of hyperbolic
components, according to the behavior of the critical points:
{\it adjacent, bitransitive, capture,} and {\it disjoint} \cite{Mil}.
We are exclusively interested in the last possibility:  
a component $\h\subset\overset{\circ}{\Cub}$ is of disjoint type $D_{m,n}$
if $F_{a,b}$ has distinct attracting periodic orbits
with periods $m$ and $n$ for every  $(a^2,b^2)\in \h$.
By definition, the $\Per_n(\lambda)-$curve consists of all
parameter values for which the cubic polynomial $F_{a,b}$
has a periodic point of period $n$ and multiplier $\lambda$.
The geography of $\Per_1(0)$ was studied in \cite{Mil-per} and
\cite{Faught}.

Notice that if the coefficients of a cubic polynomial are real
then so are the corresponding parameters $A$ and $B$. Thus we
may consider the connectedness locus of real cubic maps,
the set of pairs $(A,B)\in \BBB R^2$ such that $J(F_{a,b})$ is connected.
This locus $\CubR$ is also compact, connected and cellular \cite{Mil}.
We refer the reader to \figref{AB-plane} which was generated
by a computer program of Milnor. The real slices of various hyperbolic
components are rendered in different shades of gray. 
Certain disjoint type components are indicated, as are
the curves $\Per_1(1)\cap \CubR$ and $\Per_2(1)\cap \CubR$.

\realfig{symmlocus}{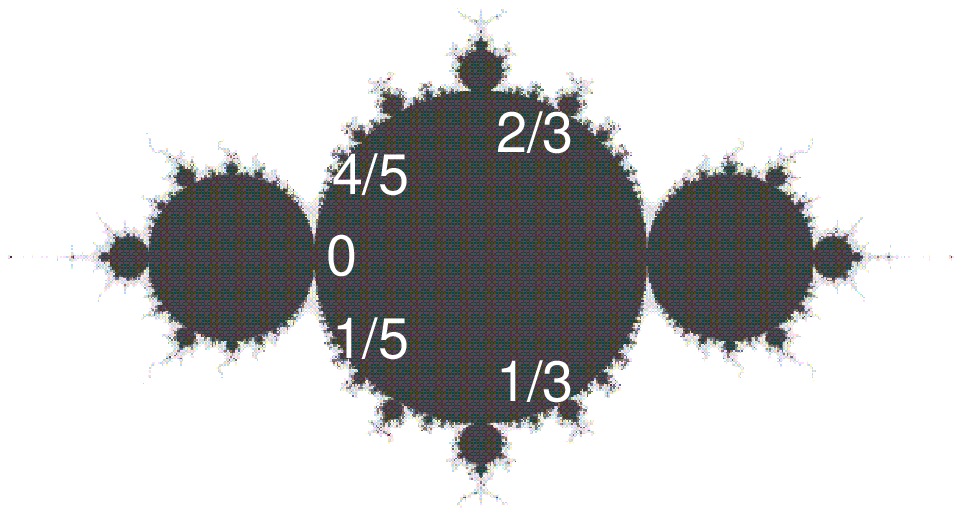}
{Symmetry locus in the family $P_{A,D}$}{10cm}

To avoid ambiguities arising from the choice of normalization,
we will actually work in the family of cubics
$$P_{A,D}=A(w^3-3w)+D, \;\;A\neq 0$$
with {\it marked critical points} $-1$ and $+1$. The
reparametrization
$$\BBB C^*\times\BBB C\ni (A,D)\mapsto (A,AD^2)=(A,B)\in \BBB C^*\times \BBB C$$
is branched over the {\it symmetry locus} $B=0$ consisting of normalized 
cubics which commute with $z\mapsto -z$ (see \figref{symmlocus}). 
In particular, 
$$\CubM=\{(A,D)\subset \BBB C^*\times \BBB C|\; 
J(P_{A,D})\;\operatorname{is\; connected}\}$$
is a branched double cover of $\Cub\cap(\BBB C^*\times \BBB C)$.
The marking of critical points allows us to 
label the attracting cycles of maps in disjoint type components
$\h\subset\CubM$, and
we denote the corresponding multipliers $\lambda^\pm_\h(A,D)$.
It is shown in \cite{Mil-components} that the maps $\Lambda_\h:
\h\to\BBB D\times\BBB D$ given by
$$\Lambda_\h(A,D)=(\lambda^-_\h(A,D)\,,\lambda^+_\h(A,D))$$
are biholomorphisms. The omitted
curve $A=0$, consisting of maps with a single degenerate critical point,
is irrelevant to the discussion of disjoint type components.

This useful change of variable has the unfortunate side-effect that the
values $(A,D)\in \BBB R^*\times \BBB R$ only account for the first
and third quadrants of the real $(A,B)-$plane, the second and fourth quadrants
being parameterized by $\BBB R^*\times i\BBB R$. We are therefore unable
to furnish a faithful illustration of the entire locus
$$\CubMR=\{(A,D)| \; (A,AD^2)\in\CubR\}.$$

\subsection{Tools}
\label{tools}
\medskip

Let $f:W'\to W$ be a quadratic-like map with connected Julia set, 
$\zeta$ a repelling fixed point with combinatorial rotation number $p/q$ and
associated quotient torus
$T_\zeta=(D\setminus\{\zeta\})/f^{\circ q}$, where $D$ is a fixed but 
otherwise arbitrary linearizing neighborhood $D\ni \zeta$. Given $S
\subset W\setminus\{\zeta\}$ with $f^{\circ q}(S\cap D)=S\cap f^{\circ q}(D)$,
we denote $\hat S$ its projection to $T_\zeta$; in particular,
$\hat{K}_1(f),\ldots,\hat{K}_q(f)\subset T_\zeta$ 
are the quotients of the various components of $K(f)\setminus \{\zeta\}$.
As any two annuli $A_1\supset \hat{K}_{i_1}(f)$ and $A_2\supset
\hat{K}_{i_2}(f)$ are isotopic we may speak of a distinguished isotopy class of
annuli $A\subset T_\zeta$, namely $A\sim \hat{K}(f)$ if and only if $A$
is isotopic to an annulus containing some $\hat{K}_i(f)$. Moreover, it
is easy to see if $A\subset \hat{K}_i(f)$ does not separate $T_\zeta$ then 
$A\sim\hat{K}(f)$; it follows then that $\zeta$ is on the boundary of an 
immediate attracting basin. Consider
$$\underline{\mod}\;\hat{K}_i(f)=\sup\{\mod A|\;A\subset \hat K_i(f)\}$$
and
$$\overline{\mod}\;\hat{K}_i(f)=\inf\{\mod A|\;A\supset \hat K_i(f)\}$$
over annuli $A\sim\hat{K}(f)$. Notice that these quantities are independent
of $i$. In view of the following we may simply write ${\mod}\;\hat{K}(f)$:

\begin{lem}
\label{equal-mod}
In this setting $\underline{\mod}\;\hat K_i(f)=\overline{\mod}\;\hat K_i(f)$.
\end{lem}

\begin{pf}
It is obvious that $\underline{\mod}\;\hat{K}\leq\overline{\mod}\;
\hat{K}$ for $\hat{K}=\hat{K}_i(f)$.  Conversely,  given $R_n\searrow 
R_\infty=e^{\overline{\mod}\; K}$
there exist conformal embeddings $h_n:{\BBB A}_{R_n}\to T_\zeta$ such that 
$h_n({\BBB A}_{R_n})\supset \hat{K}_i(f)$, where
$${\BBB A}_R=\{z: 1<|z|<R\}.$$
It follows from standard estimates in geometric function theory that
the $h_n$ form a normal family on ${\BBB A}_{R_\infty}$; moreover, 
as all of these embeddings are isotopic, every
limit $h_\infty=\lim_{k\to\infty} h_{n_k}$ is univalent. Clearly, 
$h_\infty({\BBB A}_{R_\infty})\subset \hat{K}$ and therefore \linebreak 
$\overline{\mod}\;\hat{K}\leq\underline{\mod}\;\hat{K}$.
\end{pf}

As $\underline{\mod}\;\hat{K}(f)$ is defined in terms of the interior
of $K(f)$, we observe: 

\begin{rem}
\label{indep-mod}
${\mod}\;\hat{K}(f)={\mod}\;\hat{K}(g)$ at corresponding fixed points of
hybrid equivalent quadratic-like maps $f$ and $g$.
\end{rem}

Let $f:W'\to W$ be a quadratic-like map with connected Julia set, and $\zeta$
a repelling fixed point with combinatorial rotation number $p/q$.
An {\it invariant sector} with vertex $\zeta$ is a simply connected domain
$S\subset W$ bounded by an arc of $\partial W$ and two additional arcs
$\gamma_1$ and $\gamma_2$ with $\gamma_j\subset f^{\circ q}(\gamma_j)$ and 
a common endpoint at $\zeta$. We write $S=\sector{\gamma_1}{\gamma_2}$
for the sector between $\gamma_1$ and $\gamma_2$ as listed in counterclockwise
order. The quotient $\hat S\subset T_\zeta$ is an open annulus whose
modulus will be referred to as the {\it opening modulus} $\mod S$
of the sector $S$. 
\realfig{sector-fig}{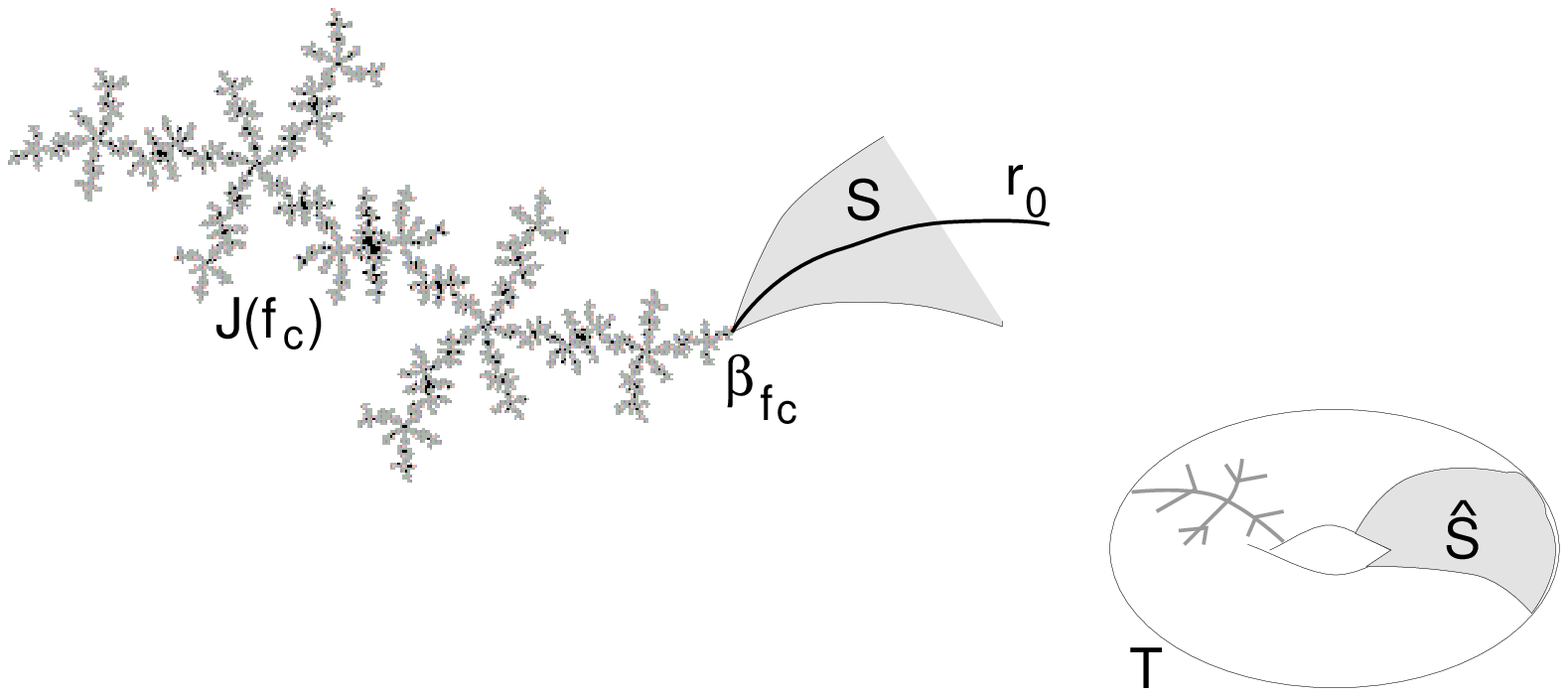}{An invariant sector $S_t(r_0)$
with vertex at $\beta_{f_c}$}{0.74\hsize}

Consider a restriction of a quadratic polynomial $f_c$ with a connected Julia
set to the domain $W\supset K(f_c)$ bounded by an equipotential $E_\rho$.
Invariant sectors for this map  may be constructed as follows
\cite{BD}: given  a ray $r_\theta$ landing at a fixed point $\zeta$ with
combinatorial rotation number $p/q$, consider
$$S_t(r_\theta)=B_{f_c}^{-1}(\{e^{r+2\pi i\gamma}|\;1<r<\rho\;\mbox{ and }\; 
|\gamma-\theta|< tr\}).$$ 
\figref{sector-fig} depicts an invariant sector $S=S_t(r_0)$ and its
quotient $\hat S$.

\begin{lem}[\cite{BF}, Prop. 4.1]
\label{disjoint-sectors}
Given $\rho>1$ there exists $\tau>0$ such that for any $t<\tau$
the domains $S_t(r_{2^{i-1}\theta})\subset W$ for $i=1,\ldots,q$ 
are disjoint invariant sectors.
\end{lem}

For the readers convenience let us review the notion of an almost
complex structure. Let $\sigma=\{E_z\}_{z\in G}$ be a measurable field 
of ellipses on a planar domain $G$ with the ratio of major to minor axes
at the point $z$ denoted by $K(z)$. 
The {\it complex dilatation} 
is a complex valued function $\mu:G\to \BBB D$, where
$|\mu(z)|=(K(z)-1)/(K(z)+1)$, and the argument of $\mu(z)$ is twice
the argument of the major axis of $E_z$.
A {\it bounded measurable almost complex structure} is a 
field of ellipses $\sigma$  with $\|\mu\|_\infty<1$.
The {\it standard almost complex structure} $\sigma_0$ is a field of
circles, thus having identically vanishing complex dilatation.

Given an ellipse field  $\sigma$ on $G$ and  
 an almost everywhere differentiable homeomorphism  $h:W\to G$
the {\it pullback of $\sigma$} is an ellipse field  $h^*\sigma$ on $W$
obtained as follows.
 For almost every $z\in W$, there is  a linear tangent
map $$T_zh:T_zW\to T_{h(z)}G.$$ Let 
$\sigma=\{E_\zeta\subset T_\zeta(G)\}_{\zeta\in G}$, then
$h^*\sigma$ is given by $\{T_zh^{-1}(E_{h(z)})\subset T_zW\}_{z\in W}$.
We note that when the map $h$ is quasiconformal
 the pullback of the standard structure $\sigma=h^*\sigma_0$
is a  bounded almost complex structure.

The proofs of the following general principles can be found in \cite{Le}:

\begin{thm} 
\label{conformal}
Let $h$ be a quasiconformal map such that $h^*\sigma_0=\sigma_0$.
Then $h$ is conformal.
\end{thm}

\begin{thm}[Measurable Riemann Mapping Theorem]
\label{MRMT}
If $\sigma$ is a bounded  almost complex structure on a domain 
$G\subset \BBB C$,  then there exists a quasiconformal homeomorphism
$h:G\to \BBB C$, such that $$\sigma=h^*\sigma_0.$$
\end{thm}

\bigskip

\section{Outline of the Results}
\label{outline}

In the picture of the real cubic connectedness locus (\figref{AB-plane})
one observes several shapes reminiscent of the Mandelbrot set 
(\figref{mandelbrot}). We quote Milnor (\cite{Mil}): 
``... these embedded copies tend to be discontinuously distorted
at one particular point, namely the period one saddle node point $c=1/4$,
also known as the root of the Mandelbrot set. The phenomenon is
particularly evident in the lower right quadrant, which exhibits
a very fat copy of the Mandelbrot set with the root point 
stretched out to cover
a substantial segment of the saddle-node curve $\Per_2(1)$.
$\ldots$
As a result of this stretching, the cubic connectedness locus
fails to be locally connected along this curve.''

The original goal of our investigation was to explain the appearance
of these distorted copies of the Mandelbrot set embedded in $\CubR$.
This has lead us to the following results:

For $p/q\in{\BBB Q}/{\BBB Z}$ we consider the set
$C_{p/q}\subset \CubM$ consisting of cubic polynomials for which
$2q$ distinct external rays with combinatorial rotation number
$p/q$ land at some fixed point $\zeta$. As there can be at most one such
point, the various $C_{p/q}$ are disjoint. 
Each $C_{p/q}$ is in turn the disjoint union of
subsets $C_{p/q,m}$ indexed by an odd integer $1\leq m\leq 2q-1$ 
specifying how many of these rays are encountered in passing counterclockwise
from the critical point $-1$ to the critical point $+1$. 
In particular, $C_0$ consists of those cubics in $\CubM$ 
whose fixed rays $r_0$ and $r_{1/2}$ land at the same fixed point.
 

\begin{thm}[Main Theorem]
\label{main}
Given $p/q$ and $m$ as above, there exists a homeomorphic embedding
$$h_{p/q,m}:\M_{p/q}\setminus\{\rt_{p/q}\}\,\times\,
\M_{p/q}\setminus\{\rt_{p/q}\}\,\longrightarrow\, C_{p/q,m}\,$$
mapping the product of hyperbolic components 
$H_{p/q}\times H_{p/q}$ onto a component $\h_{p/q,m}$ of type $D_{q,q}$.
\end{thm}

We note that $\h_{p/q,m}$ is the unique $D_{q,q}$ component contained in
$C_{p/q,m}$ as will follow from \thmref{renorm-inj}.
The restriction of $h_{p/q,m}$ to $H_{p/q}\times H_{p/q}$ is easily
expressed in terms of the multiplier maps defined in \secref{preliminaries}:
$$h_{p/q,m}(\rho,\tl\rho)=\Lambda_{\h_{p/q,m}}^{-1}(\,\lambda_{H_{p/q}}
(\rho),\lambda_{H_{p/q}}(\tl \rho)\,).$$ Discontinuity of $h_{p/q,m}$ at
the corner point $(1,1)$ is a special case of a phenomenon studied by one
of the authors: 

\begin{thm}\cite{mate}
\label{two-disks}
Each algebraic homeomorphism 
$$\Lambda_{\h_{p/q,m}}:\h_{p/q,m}\to \BBB D\times \BBB D$$
extends to a continuous surjection
$\bar{\h}_{p/q,m}\to\bar{\BBB D}\times\bar{\BBB D}$.
The fiber over $(1,1)$
is the union of two closed discs whose boundaries are real-algebraic curves 
with a single point in common, and all other fibers are points.
\end{thm}

The following reasonable conjecture appears to be inaccessible by purely
quasiconformal techniques:

\begin{conjecture} \label{wishful-thinking}
Each $h_{p/q,m}$ extends to a continuous embedding
$$\M_{p/q}\times\M_{p/q}\,\setminus\,
\{(\rt_{p/q},\rt_{p/q})\}\,\longrightarrow\,C_{p/q,m}\,.$$
\end{conjecture}

\medskip
We draw additional conclusions from the natural symmetries of our construction.
The central disk in \figref{symmlocus} is parameterized by the
eigenvalue $-3A$ of the attracting fixed point
at $0$; this region corresponds to symmetric cubics whose Julia sets are
quasicircles. Each value $A=-\frac{1}{3}e^{2\pi ip/q}$ yields a
map with a parabolic fixed point at 0. These parameters
are evidently the roots of small embedded copies of $\M$, and our results
confirm this observation for odd-denominator rationals.
More specifically, it will follow that the latter copies are the
images of $h_0\circ\Delta$ and $h_{p/q,q}\circ\Delta$
for odd $q>1$, where $$\BBB C \ni c\overset{\Delta}{\mapsto}(c,c)\in\BBB C^2$$
is the diagonal embedding. As $P^{\circ 2}_{A,0}$ and $P^{\circ 2}_{-A,0}$ are
conjugate by $z\mapsto -z$, our construction also accounts for the copies
with $q\equiv 2$ (mod $4)$.
Every map in the symmetry locus is semiconjugate, via the quotient
determined by the involution, to a cubic polynomial with a fixed critical
value. Such maps were studied by Branner and Douady \cite{BD} who effectively
prove that the entire limb attached at the parameter value $A=-1/3$ is
a homeomorphic copy of the limb $L_{1/2}\subset \M$; it can be shown by a
variant of the pullback argument in \secref{renormalization}
that the image of $h_0\circ
\Delta$ corresponds to the small copy $\M_{1/2}\subset L_{1/2}$.  

Similar considerations applied to the antidiagonal embedding yield results
for the real connectedness locus.
In view of the fact that real polynomials commute with 
complex conjugation, $\CubMR\cap C_{p/q}=\emptyset$ unless $p/q\equiv-p/q\;
\mbox{(mod 1)}$, and it therefore suffices to consider the real slices of 
$C_{1/2}$ and $C_0$.

\begin{thm}
\label{real-main-thm}
There exist homeomorphic embeddings
$$\Psi_{1/2,1}:\M_{1/2}\setminus\{\rt_{1/2}\}\to
\CubMR\cap C_{1/2,1}$$
$$\Psi_{1/2,3}:\M_{1/2}\setminus\{\rt_{1/2}\}\to
\CubMR\cap C_{1/2,3}$$
and
$$\Psi_0:\M\setminus\{\rt\}\to\CubMR\cap C_0.$$
\end{thm}

\medskip
It follows from recent work of Buff \cite{Buff} 
that these maps are compatible with the standard embeddings in 
the plane (see the discussion in \secref{conclusions}).
Their projections in $\CubR$ are indicated in \figref{AB-plane}.
Notice that the two images of $\M_{1/2}\setminus\{\rt_{1/2}\}$
have been identified
while the image of $\M\setminus\{\rt\}$ has been folded in half.
The latter defect
is overcome through passage to the $(A,\sqrt{B})$-plane, at the cost of 
both copies of $M_{1/2}\setminus\{\rt_{1/2}\}$; we thank
J. Milnor for enabling us to include \figref{per11comb} where the 
comb on the $D_{1,1}$ component is better resolved. The existence
of this comb is verified with the aid of techniques developed by Lavaurs
\cite{Lav}.

\begin{thm}
\label{non-loc-conn}
The real cubic connectedness locus is not locally connected.
\end{thm}

The remainder of this paper is structured as follows.
In \secref{intertwining} we construct cubic polynomials by means of
{\it quasiconformal surgery} on pairs of quadratics. The issues of
uniqueness and continuity are addressed in \secref{renormalization}
through the use of the renormalization operators $\lren$ and $\rren$ 
defined for {\it birenormalizable} cubics; together they essentially
invert the surgery. We show that $\lrren$ is a homeomorphism and then
complete the proofs of Theorems \ref{main} and \ref{real-main-thm} in
\secref{conclusions}.
The measurable dynamics of birenormalizable maps
is discussed in  \secref{res-mes}. In Appendix \ref{discont-root} we
comment further on the discontinuity described in  
\thmref{two-disks}, and we conclude by proving \thmref{non-loc-conn}
in Appendix \ref{section-comb}.

It is worth noting that quasiconformal surgery is only employed in
the proof of surjectivity for $\lrren$.
More generally, we might associate a pair of renormalization operators to 
any disjoint type hyperbolic component $\h\subset \CubM$ in the hope
of finding an embedded product of Mandelbrot sets ``growing'' from $\h$,
but we are  unable to adapt our surgery construction to this
broader setting. Part II of this paper
will present a different approach to proving surjectivity of birenormalization,
culminating in a more general version of \thmref{main}.

\section{Intertwining surgery}
\label{intertwining}
\subsection{History} 
The intertwining construction was described
in the 1990 Conformal Dynamics Problem List \cite{Biel}:
``Let $P_1$ be a monic polynomial with connected Julia set having a repelling
fixed point $x_0$ which has ray landing on it with rotation number
$p/q$. Look at a cycle of $q$ rays which are the forward images of the first.
Cut along these rays and get $q$ disjoint wedges. Now let $P_2$
be a monic polynomial with a ray of the same rotation number landing 
on a repelling periodic point of some period dividing $q$ 
(such as $1$ or $q$). Slit this dynamical plane along the same rays making
holes for the wedges. Fill the holes in by the corresponding wedges above 
making a new sphere. The new map is given by $P_1$ and $P_2$, 
except on a neighborhood of the inverse images of the cut rays where it
will have to be adjusted to make it continuous.''

\subsection{Construction of a cubic polynomial}
Fix $p/q$ written in lowest terms and an odd integer $m=2k+1$ between 1 and
$2q-1$. Our aim is to construct a map 
$$h_{p/q,m}:\M_{p/q}\setminus\{\rt_{p/q}\}\,\times\,
\M_{p/q}\setminus\{\rt_{p/q}\}\longrightarrow C_{p/q,m}.$$

Fixing parameter values $c$ and $\tl c$ in $\M_{p/q}\setminus\{\rt_{p/q}\}$,
consider quadratic-like maps  $f:W'\to W$ and $\tl f:\tl W'\to \tl W$
hybrid equivalent to $f_c$ and $\tl f_c$ respectively, the choice
of the hybrid equivalences to be made below.
In what follows we will identify $W$ and $\tl W$
to obtain a new surface. The reader is invited
to follow the construction in the particular case $p/q=1/2$ with
$m=3$, as illustrated in \figref{int-figure}.

We lose no generality by
assuming that $0$ is the critical point for both $f$ and $\tl f$.
Let $\zeta$ be the unique repelling fixed point of $f$ with combinatorial
rotation number $p/q$, that is $\zeta=\beta_f$ for $p/q=0$ and 
$\zeta=\alpha_f$ otherwise, and
$S_i\equiv\sector{l_i}{r_i},\;i=0,\ldots,q$ in $W\setminus K(f)$ 
a cycle of disjoint invariant sectors 
with vertex $\zeta$, indexed in counterclockwise order so that 
the critical point $0$ lies in the complementary region
between $S_{q-1}$ and $S_0$. We similarly specify $\tl \zeta$
and a cycle of invariant sectors $\tl S_i$ for $\tl f$. Let 
$$\phi:\bigcup_{i=1}^q l_i\cup r_i\to \bigcup_{i=1}^q \tl l_i\cup \tl r_i$$
sending $l_i$ to $\tl r_{i+q-k}$ and $r_i$ to
$\tl l_{i+q-k-1}$, where $k=(m-1)/2$ and indices are understood modulo $q$,
be any smooth conjugacy,
$$\phi(f(z))={\tl f}(\phi(z)).$$

\realfig{int-figure}{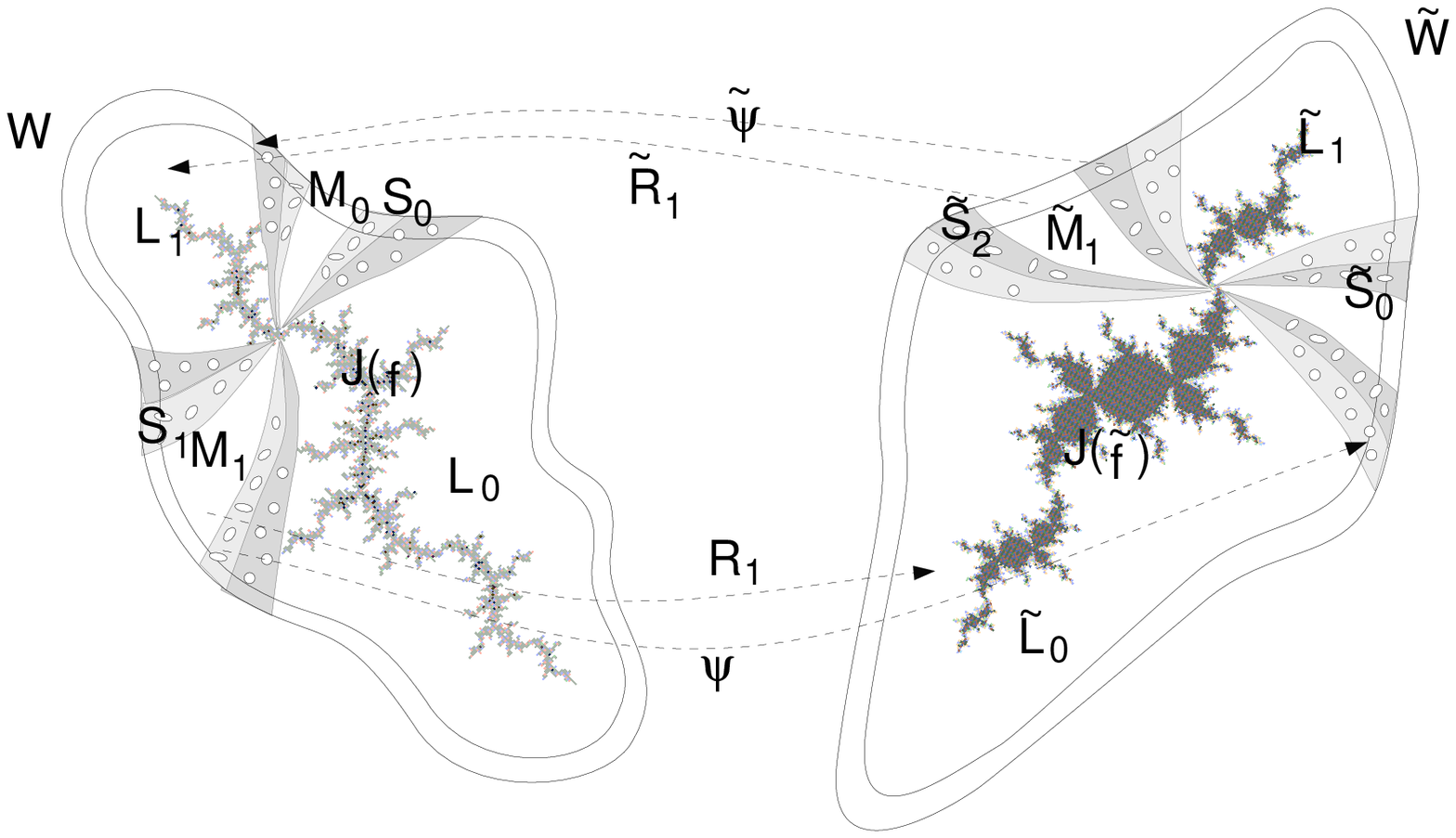}{Superimposing the regions $W$ and
$\tl W$}{13cm}

The sector $S_i$ should now correspond to the component of 
$K(\tl f)\setminus\{\tl\zeta\}$
containing $\tl f^{\circ i+q-k}(0)$. An informal rule
known as {\it Shishikura's Principle} warns against altering the
conformal structure on regions of uncontrolled recurrence, and we will 
therefore employ invariant sectors $M_i\subset S_i\subset N_i\subset 
W\setminus K(f)$ and $\tl M_i\subset\tl S_i\subset\tl M_i\subset
\tl W\setminus K(\tl f)$ to be determined below. For $i=0,\ldots,q$ denote
$L_i$ the component of $W\setminus \bigcup_{j=1}^q N_j$ 
containing $f^{\circ i}(0)$, and $\tl L_i$ the corresponding 
component of $\tl W\setminus \bigcup_{j=1}^q \tl N_j$.
Let $$R_i:M_i\to\tl L_{i+q-k}$$ be the conformal homeomorphism
sending $\partial M_i\cap W$ to 
$\partial \tl L_{i+q-k}\cap\tl W$ so that $\zeta$ maps to $\tl \zeta$, and
$$\tl R_i:\tl M_i\to L_{i-q+k+1}$$ the conformal homeomorphism 
sending $\partial \tl M_i\cap \tl W$ to 
$\partial L_{i-q+k+1}\cap W$ so that $\tl\zeta$ maps to $\zeta$.
These Riemann maps
extend continuously to the sector boundaries, and it remains to fill in
the gaps:

\begin{prop}[Quasiconformal interpolation]
\label{interpolation}
For any pair  $f_c$ and $f_{\tl c}$ as above  the hybrid equivalent
quadratic-like maps $f$ and $\tl f$ and the invariant
sectors $$N_i\supset S_i\supset M_i\;\mbox{ and }\;
\tl N_i\supset\tl S_i\supset \tl M_i$$
may be chosen so that there exist quasiconformal maps
$$\psi:\bigcup_{i=1}^q (S_i\setminus M_i)\to \bigcup_{i=1}^q (\tl N_i\setminus
\tl S_i)$$
$$\tl \psi:\bigcup_{i=1}^q (\tl S_i\setminus \tl M_i)\to \bigcup_{i=1}^q
( N_i\setminus  S_i)$$
with 
$$\psi|_{\partial M_i}=R_i|_{\partial M_i}\;\mbox{ and }\;
\psi|_{\partial S_i}=\phi|_{\partial S_i}$$
$$\tl \psi|_{\partial\tl M_i}=\tl R_i|_{\partial \tl M_i}\;\mbox{ and }\;
\tl\psi|_{\partial\tl S_i}=\tl \phi|_{\partial \tl S_i}.$$
\end{prop}

\medskip
Let us complete the construction assuming the truth of 
Proposition \ref{interpolation}.
We choose sectors $N_i$, $M_i$, $\tl N_i$, $\tl M_i$,
and maps $\psi$, $\tl \psi$ as specified above.
These identifications turn $W\sqcup \tl W$ into a new
manifold with a natural quasiconformal atlas.
Consider the almost complex structure
$\sigma$ on $W$, given by $\psi^*\sigma_0$ on $\bigcup_{i=1}^q S_i
\setminus M_i$ and by $\sigma_0$ elsewhere; similarly, let $\tl \sigma$ be the 
almost complex structure on $\tl W$ given by
$\tl \psi^*\sigma_0$ on $\bigcup_{i=1}^q \tl S_i\setminus \tl M_i$,
and $\sigma_0$ 
elsewhere. In view of \thmref{MRMT} there exist quasiconformal homeomorphisms
$h:W\to X$ and $\tl h:\tl W\to \tl X$ such that
$\sigma=h^*\sigma_0$ and $\tl \sigma=\tl h^*\sigma_0$.
The domains $X$ and $\tl X$ give coordinate neighborhoods,
and the maps $\psi\circ h^{-1}$, $R_i\circ h^{-1}$, $\tl \psi\circ\tl  h^{-1}$,
and $\tl R_i\circ \tl h^{-1}$ yield an atlas of a Riemann surface
with the conformal type of a punctured disc. We obtain a conformal disc
$\Delta$ by replacing the puncture with a point $\star$. Setting
$$\Delta'=h(W'\setminus \bigcup_{i=1}^q S_i)\cup \tl h(\tl W'\setminus 
\bigcup_{i=1}^q \tl S_i)\cup\{\star\},$$ 
we define a new map $F:\Delta'\to \Delta$ by 
$$F(z)=\left\{ \begin{array}{lll}
        h(f(h^{-1}(z))) & \mbox{for} &
           z\in h(W'\setminus \bigcup_{i=1}^q S_i),\\
 	\tl h(\tl f(\tl h^{-1}(z))) & \mbox{for} &
           z\in \tl h(\tl W'\setminus \bigcup_{i=1}^q \tl S_i),\\
        \star & \mbox{for} & z\in\{\star,-\zeta,-\tl\zeta\}.
	\end{array}\right.$$
It is easily verified that $F$ is a three-fold branched covering with 
critical points $h(0)$ and $\tl h(0)$, and analytic except on the
preimage of
$${\cal S}=\bigcup_{i=1}^q h(S_i\setminus M_i)\cup\tl h(\tl S_i\setminus
\tl M_i).$$

Recalling that the sectors $S_i$ and $\tl S_i$ are invariant
and disjoint, we consider the following almost complex structure
on $\Delta$:
$$\left\{\begin{array}{ll}
		\hat{\sigma}=(F^{\circ n})^*\sigma_0 & \mbox{on }
   		F^{\circ -n}({\cal S})\\ 
		\hat{\sigma}=\sigma_0 & \mbox{elsewhere }.
		\end{array}\right.$$ 
By construction, the complex dilatation of $\hat{\sigma}$ has the same bound as
that of $F^*\sigma_0$, and moreover
$$F^*\hat{\sigma}=\hat{\sigma}.$$
It follows from \thmref{MRMT} that there is a quasiconformal homeomorphism
\linebreak 
$\varphi:\Delta\to V\subset\BBB C$ with $\hat{\sigma}=\varphi^*\sigma_0$.
Setting $U=\varphi(\Delta')$, we obtain a cubic-like map
$$G=\varphi\circ F\circ \varphi^{-1}:U\to V.$$

In view of \thmref{strt-thm} there is a unique hybrid equivalent
cubic polynomial $P_{A,D}$ whose critical points $-1$ and $+1$ correspond
to the critical point of $f$ and $\tl f$ respectively.
The construction yields extensions of the natural embeddings
$$\pi:K(f_c)\to K(P_{A,D})\;\mbox{ and }\;
\tl\pi:K(f_{\tl c})\to K(P_{A,D})$$
to neighborhoods of the filled Julia sets. 

\begin{rem}
\label{conf-julia}
By construction, the projections $\pi$ and $\tl \pi$ are 
conformal on the respective filled Julia sets:
$$\bar\partial\pi=0 \mbox{ a.e. on } K(f_c)\;\mbox{ and }\;
\bar\partial{\tl \pi}=0 \mbox{ a.e. on } K( f_{\tl c}).$$
\end{rem}

\medskip
We write $P_{A,D}\approx f_c{\underset{p/q,m}{\curlyvee}}f_{\tl c}$ for 
{\it any} cubic polynomial so obtained. It is not yet clear that this 
correspondence is well-defined, let alone continuous. 
These issues will be addressed in \secref{renormalization}.

\subsection{Quasiconformal interpolation}
{\bf Proof of Proposition \ref{interpolation}:}\\
Note first that were it not for the condition of quasiconformality,
the existence of the interpolating maps $\psi$ and $\tl \psi$ would 
follow without any additional argument. Any smooth
interpolations are quasiconformal away from the points 
of $\zeta$ and $\tl \zeta$, the issue is the compatibility of the 
the local behaviour of $\phi$ at $\zeta$  with that of 
$R_i$ and $\tl R_i$.

\begin{lem}
\label{choice-sectors}
Given any $c\in\M_{p/q}\setminus\{\rt_{p/q}\}$ and $\upsilon>0$, there exists
a quadratic-like map $f$ which is hybrid equivalent to $f_c$ and 
admits disjoint invariant sectors $S_i$ as above with ${\mod}\; S_i>\upsilon$.
\end{lem}
 
\begin{pf}
We begin by fixing a quadratic-like restriction $f_c: G'\to G$ between
equipotentially bounded regions, and apply \lemref{disjoint-sectors} to
obtain a cycle of disjoint invariant sectors $S_t(r_i)\subset G$.
Let $\varphi$ be a quasiconformal homeomorphism from the annulus
$\hat S_t(r_0)\subset T_\zeta$ to some standard annulus ${\BBB A}_\rho$
with $\rho>e^\upsilon$.
The almost complex structure $\sigma=\varphi^*\sigma_0$ on
$\hat S_t(r_0)$ lifts to an almost complex structure on the sector $S_t(r_0)$.
We extend this structure by pullback to the various $S_t(r_i)$ and their
preimages, and extend by $\sigma_0$ elsewhere, to obtain an invariant almost 
complex structure $\bar\sigma$ on $G$. In view of \thmref{MRMT} there
exists a quasiconformal homeomorphism $\varphi:G\to \varphi(G)\subset{\BBB C}$
with $\bar\sigma=\varphi^*\sigma_0$, giving a hybrid equivalence between
$f_c$ and the quadratic-like map
$$f=\varphi\circ f_c\circ \varphi^{-1}:\varphi(G')\to \varphi(G).$$
It follows from \thmref{conformal} that ${\mod}\; S_i=\upsilon$ where $$S_i=
\varphi(S_t(r_i)).$$
\end{pf}

Given $c,\tl c\in\M_{p/q}\setminus\{\rt_{p/q}\}$, we apply
\lemref{choice-sectors} to $f_c$ and $f_{\tl c}$ to obtain hybrid equivalent
quadratic-like maps $f$ and $\tl f$ admitting invariant sectors $S_i$ and
$\tl S_i$ with  
$${\mod}\;S_i >{\mod}\;\hat{K}(f_{\tl c})\;\mbox{ and }\;
{\mod}\;\tl S_i >{\mod}\;\hat{K}(f_c).$$
In view of Remark \ref{indep-mod} we may then choose disjoint 
invariant sectors $N_i\supset S_i$ and $\tl N_i\supset \tl S_i$
so that 
$${\mod}\; S_i >{\mod}\;\tl L_j\;\mbox{ and }\; {\mod}\; \tl S_i >{\mod}\;L_j$$
for the complementary invariant sectors $L_j$,$\tl L_j$  as above. Finally,
we choose $M_i\subset S_i$ and $\tl M_i\subset \tl S_i$ with 
$${\mod}\;M_i={\mod}\; \tl L_j\;\mbox{ and }\;{\mod}\;\tl M_i={\mod}\;L_j.$$

We now exploit the following observation of \cite{BD}; see
\cite[Lemmas 6.4, 6.5]{Biel-thesis} for a detailed exposition.

\begin{lem}
With this choice of maps and invariant sectors 
there exist desired quasiconformal interpolations

$$\psi:\bigcup_{i=1}^q (S_i\setminus M_i)\to \bigcup_{i=1}^q (\tl N_i\setminus
\tl S_i)$$
$$\tl \psi:\bigcup_{i=1}^q (\tl S_i\setminus \tl M_i)\to \bigcup_{i=1}^q
( N_i\setminus  S_i)$$
with
$$\psi|_{\partial M_i}=R_i|_{\partial M_i}\;\mbox{ and }\;
\psi|_{\partial S_i}=\phi|_{\partial S_i}$$
$$\tl \psi|_{\partial\tl M_i}=\tl R_i|_{\partial \tl M_i}\;\mbox{ and }\;
\tl\psi|_{\partial\tl S_i}=\tl \phi|_{\partial \tl S_i}.$$
\end{lem}

\section{Renormalization}
\label{renormalization}

\subsection{Birenormalizable cubics}
Throughout this section we will work with fixed values of $p/q$ and $m$
as specified above. Here we describe the construction which will 
provide the inverse to the map $h_{p/q,m}$.

We start with a cubic polynomial $P=P_{A,D}$ with
$(A,D)\in C_{p/q,m}$. Let $\zeta$ be the landing point of 
the periodic rays with rotation number $p/q$, and denote
$K_{\pm1}$ the connected components of $K(P)\setminus \{\zeta\}$
containing the critical points $\pm 1$.
Below we determine  quadratic-like restrictions of $P$  which
domains contain the appropriate critical points. To fix the
ideas, we illustrate this {\it thickening procedure} for left renormalizations 
only.

Let $r_{\theta_1}$ and $r_{\theta_2}$ be the two periodic external 
rays landing at $\zeta$ which separate $K_{-1}$ from the other rays landing
there; without loss of generality, $0\leq\theta_2<\theta_1<2\pi$ so that 
the rays landing at $K_{-1}$ have angles in $[\theta_2,\theta_1]$.
Choose a neighborhood $U\ni\zeta$ corresponding to
a round disc in the local linearizing coordinate. Fix an equipotential
$E$ and a small $\eps>0$, and consider the segments of the 
rays $r_{\theta_1+\eps}$ and $r_{\theta_2-\eps}$ connecting the
boundary of $U$ to $E$. Let $\Omega\supset K_{-1}$ be the region bounded
by these two ray segments and the subtended arcs of $E$ and $\partial U$,
and consider the component $\Omega'$ of $P^{\circ -q}(\Omega)$ with 
$\Omega'\subset\Omega$. In view of the fact that $\zeta$ is repelling,
$\bar\Omega'\subset\Omega$ provided that $\eps$ is sufficiently small.
Thus, $$P^{\circ q}:\Omega'\to \Omega$$
is a quadratic-like map which filled Julia set will be denoted $K_\lren$.
This set is connected  if and only if 
$\{P^{\circ nq}(-1)\}_{n=0}^\infty\subset K_{-1}$, in which case we
 refer to the unique hybrid conjugate 
quadratic polynomial $f_c$ as the {\it left renormalization} $\lren(P)$
and  call   $P$ {\it renormalizable to the left}.

\figref{ren-fig} illustrates this construction for a cubic polynomial in $C_0$.
Notice that $\zeta$ becomes the $\beta-$fixed point of the new quadratic
polynomial. 
The polynomial $P$ is {\it renormalizable to the right} if 
$\{P^{\circ nq}(+1)\}_{n=0}^\infty\subset K_{+1}$, and the  set $K_{\rren}$ 
and the {\it right renormalization} $\rren(P)$ are correspondingly defined.
It follows from general considerations discussed in \cite{McM} that the
left and right renormalizations do not depend on the choice of thickened
domains. 

\realfig{ren-fig}{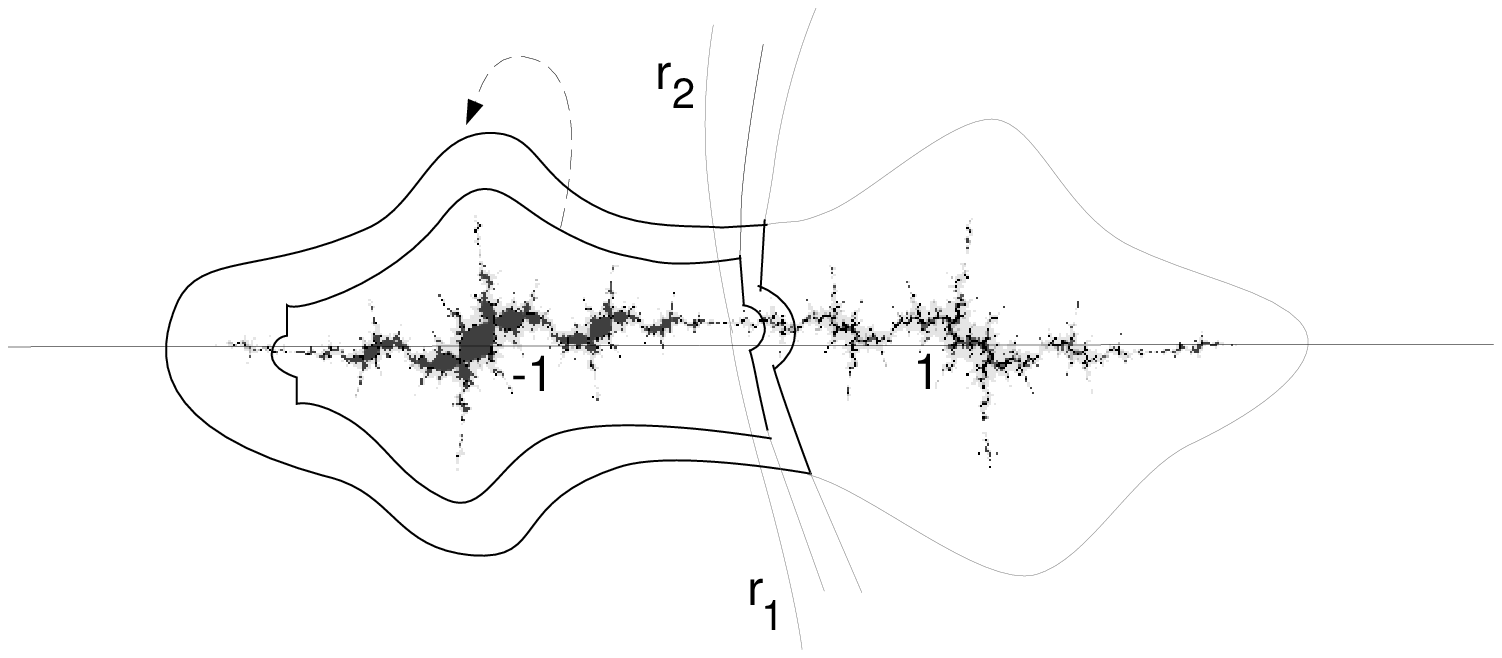}{Construction of the left
renormalization for a cubic in $C_0$}{0.8\hsize}

A cubic polynomial $P$ is said to be {\it birenormalizable} if it is
renormalizable on both left and right, in which case
\begin{equation} \label{bi-omega} 
\omega_P(-1)\cap \omega_P(+1)\subset K_{\lren}\cap K_{\rren}=\{\zeta\}
\end{equation}
and we set
$$\check{K}(P)=\bigcup_{i=0}^{\infty}P^{-\circ i}(K_\lren\cup K_{\rren}).$$
The following is an easy consequence of Kiwi's Separation \thmref{separation}
and the standard classification of Fatou components:

\begin{lem} \label{check-dense}
Let $P$ be a birenormalizable cubic polynomial. Then $\check{K}(P)$ is
dense in $K(P)$, and every periodic orbit in $K(P)\setminus\check{K}(P)$
is repelling.
\end{lem}

We denote $B_{p/q,m}$ the set of
birenormalizable cubics in $C_{p/q,m}$, writing 
$$\lrren:B_{p/q,m}\to(\M\setminus\{\rt\})\times(M\setminus\{\rt\})$$
for the map $(A,D)\mapsto(c,\tl c)$ where $f_c=\lren(P_{A,D})$ and $f_{\tl c}
=\rren(P_{A,D})$. In view of \lemref{same-rays} the thickening construction 
may be performed so that the domains of the left
and right quadratic-like restrictions vary continuously for $(A,D)\in 
B_{p/q,m}$. Applying \thmref{straightening} we obtain:

\begin{prop}
\label{renorm-cont} 
$\lrren:B_{p/q,m}\to(\M\setminus\{\rt\})\times(\M\setminus\{\rt\})$
is continuous.
\end{prop}

\medskip
The significance of intertwining rests in the following:

\begin{prop}
\label{inverse}
$\lrren:B_{p/q,m}\to(\M\setminus\{\rt\})\times(\M\setminus\{\rt\})$
is surjective.
\end{prop}

\begin{pf}
Fixing $b,\tl b\in \M\setminus\{\rt\}$, let $c=\kappa_{p/q}^{-1}(b)$ and
$\tl c=\kappa_{p/q}^{-1}(b)$ where $$\kappa_{p/q}:\M_{p/q}\to \M$$ is the
homeomorphism described in \secref{quadratics}. We saw above that
$$f_c{\underset{p/q,m}{\intr}}f_{\tl c}\approx P$$ for some cubic polynomial
$P=P_{A,D}$, and we show here that $\lrren(A,D)=(b,\tl b)$; more precisely,
we prove that $\lren(P)=f_b$, the argument for right renormalization 
being completely parallel.

Let $K(f_c)\subset W{\overset{\pi}{\rightarrow}}\BBB C$ be as in
\secref{intertwining}. The standard thickening construction yields a
quadratic-like restriction $f_c^{\circ q}:G'\to G$ with
connected filled Julia set $K\subset G\subset W$; as $f_c^{\circ q}|_{G'}
\he f_b$ it suffices to show that
this quadratic-like map is hybrid equivalent to 
$P^{\circ q}|_{\Omega'}$. By construction, $\pi$ is a quasiconformal
map conjugating $f_c^{\circ q}|_K$ to $P^{\circ q}|_{K_{\lren}}$ and 
$\bar\partial\pi(z)=0$ for almost every $z\in K$.
Let $\varphi_0:G\to\Omega$ be a quasiconformal homeomorphism with
$$\varphi_0\circ f_c^{\circ q}|_{\partial G'}=P^{\circ q}\circ \varphi_0|_
{\partial\Omega'}$$
which agrees with $\pi$ on a small neighborhood of $K$.
As $\varphi_0$ maps the critical value of $f_c^{\circ q}|_{\partial G'}$
to the critical value of $P^{\circ q}|_{\partial\Omega'}$
there is a unique lift $\varphi_1:G'\to\Omega'$ such that
$$\begin{CD}
G' @>\varphi_1>> \Omega'\\
@V{f_c^{\circ q}}VV  @V{P^{\circ q}}VV\\
G @>q_0>> \Omega
\end{CD}$$
commutes and
$\varphi_1|_{\partial G'}=
\varphi_0|_{\partial G'}$. Setting $\varphi_1(z)=\varphi_0(z)$
for $z\in \Omega\setminus \Omega'$, we obtain a quasiconformal homeomorphism
$\varphi_1:G\to\Omega$ with the same dilatation bound as $\varphi_0$; moreover,
$\varphi_1|_K=\pi|_K$. Iteration of this procedure yields a
a sequence of quasiconformal homeomorphisms $\varphi_n:G\to\Omega$ with 
uniformly bounded dilatation. The $\varphi_n$ stabilize pointwise on $G$,
 so there is a limiting quasiconformal homeomorphism
$\varphi:G\to \Omega$. By construction, 
$$\varphi\circ f_c^{\circ q}|_{G'}=P^{\circ q}\circ \varphi|_{\Omega}$$
and furthermore $\varphi|_K=\pi|_K$; it follows from Bers'
\lemref{rickman} that $\varphi$ is a hybrid equivalence.
\end{pf}

\subsection{Properness} 

Here we deduce the properness of birenormalization from Kiwi's
Separation \thmref{separation}. 

\begin{prop} \label{proper}
$\lrren:B_{p/q,m}\to(\M\setminus\{\rt\})\times(\M\setminus\{\rt\})$
is proper.
\end{prop}

\smallskip
In view of the compactness of the connectedness loci, it
suffices to prove that if $(A_k,D_k)\in B_{p/q,m}$ with
$$(A_k,D_k)\to(A_\infty,D_\infty)\in \CubM\;\;\mbox{ and }\;\;
\lrren(A_k,D_k)\to(c_\infty,\tl c_\infty)\in \M\times\M$$
then $(A_\infty,D_\infty)
\in B_{p/q,m}$ if and only if $c_\infty\neq\rt\neq\tl c_\infty$.
We require an auxiliary lemma and some further notation.
Let $\zeta_k$ be the unique repelling fixed point of $P_k=P_{A_k,D_k}$
where $2q$ external rays land, and let $\zeta_k^-$ and $\zeta_k^+$ be
the points of period $q$ which renormalize to $\alpha_{f_{c_k}}$
and $\alpha_{f_{\tl c_k}}$. We write $\mu_k$ and $\mu^\pm_k$ for the
multipliers of $\zeta_k$ and $\zeta^\pm_k$, and $\lambda_k^\pm$ for the
multipliers of the corresponding $\alpha$-fixed points. 
Passing to a subsequence if necessary, we may assume without loss of
generality that the $\zeta_k$ converge to a fixed point $\zeta_\infty$
of $P_\infty=P_{A_\infty,D_\infty}$ with multiplier $\mu_\infty$,
and the $\zeta_k^\pm$ converge to periodic points $\zeta^\pm_\infty$ with
multipliers $\mu^\pm_\infty$.

\begin{lem}
\label{stays-repelling}
In this setting, if $c_\infty\neq\rt\neq\tl c_\infty$ then 
$\zeta^{-}_\infty$, $\zeta_\infty$ and $\zeta^{+}_\infty$ belong to 
disjoint orbits. Moreover, the fixed point $\zeta_\infty$ is repelling.
\end{lem}

\begin{pf}
It follows from the Implicit Function Theorem that these orbits are
distinct unless one of $\zeta^\pm_\infty$ is parabolic with multiplier 1.
Without loss of generality $\mu_\infty^-=1$, and we may further assume 
that either $|\mu^-_k|\leq 1$ for every $k$ or $|\mu^-_k|>1$ for every
$k$. In the first case, $\lambda^-_k=\mu^-_k$ by Remark \ref{eigrmk}, 
and $\lambda^-_k\to 1$ implies $c=\rt$. In the second case, it similarly 
follows that $|\lambda^-_k|>1$ for every $k$; in view of Yoccoz Inequality, 
$c_k\in L_{p_k/q_k}$ where $p_k/q_k\to 0$, whence
$\lambda^-_k\to 1$ and again $c=\rt$. 

Because $\zeta_\infty^-$, $\zeta_\infty$ and $\zeta^\pm_
\infty$ lie in distinct orbits, it follows from \thmref{separation} 
that at least one of these points is repelling.
Suppose first that $\zeta^-_\infty$ is repelling, and let $\gamma^-_k
\in K(P_k)$ be the points which renormalize to
$-\alpha_{f_{c_k}}\in K(f_{c_k})$. Then $\gamma^-_k\to\gamma^-_\infty$
where $P_\infty(\gamma^-_\infty)=\zeta^-_\infty$, and 
the rays landing at $\gamma^-_\infty$ separate $\zeta_\infty$ from the
critical point $-1$. Similarly, if $\zeta^+_\infty$
is repelling then the rays landing at the corresponding point
$\gamma^+_\infty$ separate $\zeta_\infty$ from $+1$. 
Applying \thmref{separation} once again, we conclude that $\zeta_\infty$ is
repelling.
\end{pf}

Continuing with the proof of \propref{proper}, we observe by
\lemref{same-rays} that $\zeta_\infty$ is the common landing point of
the same two cycles of rays with rotation number $p/q$.
The thickening procedure yields a pair of quadratic-like restrictions
$$P_\infty^{\circ q}:\Omega_\infty'\to \Omega_\infty\;\;
\n\;\; P_\infty^{\circ q}:\tl\Omega_\infty'\to\tl\Omega_\infty,$$
and we may arrange for $\Omega'_\infty$ and $\tl\Omega'
_\infty$ to be the limits of thickened domains $\Omega'_k$ and $\tl\Omega'_k$
for the quadratic-like restrictions of $P^{\circ q}_k$. As 
$P_k^{\circ qn}(-1)\in \Omega'_k$ and $P_k^{\circ qn}(+1)\in \tl\Omega'_k$,
it follows that $P_\infty^{\circ qn}(-1)\in\Omega'_\infty$ and $P_\infty^
{\circ qn}(+1)\in\tl\Omega'_\infty$. Thus, $P_\infty$ is birenormalizable,
that is, $(A_\infty,D_\infty)\in B_{p/q,m}$. $\;\;\Box$

\subsection{Injectivity}

The time has come to show that the intertwining operations 
$$(f,\tl f)\mapsto f{\underset{p/q,m}{\intr}}\tl f$$ are well-defined:

\begin{prop}
\label{renorm-inj}
$\lrren:B_{p/q,m}\to(\M\setminus\{\rt\})\times(\M\setminus\{\rt\})$
is injective.
\end{prop}

\medskip
The relevant {\it pullback argument} is formalized as:

\begin{lem}
\label{pull-back}
Let $P=P_{A,D}$ and $Q=P_{\tl A,\tl D}$ where
$(A,D),\,(\tl A,\tl D)\in B_{p/q,m}$. If $$\lrren(A,D)=\lrren(\tl A,\tl D)$$
then there exists a quasiconformal homeomorphism $\varphi:\BBB C\to \BBB C$
conjugating $P$ to $Q$ with $\bar\partial \varphi=0$
almost everywhere on $\check{K}(P)$ 
\end{lem}
\begin{pf}
We begin by once again restricting $P$ and $Q$ to domains $D\supset K(P)$
and $\tl D\supset K(Q)$ bounded by equipotentials. Our first goal is the
construction of a quasiconformal homeomorphism $\varphi_0$
which is  illustrated in \figref{pull-back-fig} for $p/q=1/2$ and $m=3$. 
Let $r_1,\ldots,r_{2q}$ be the rays landing at $\zeta$, enumerated  
in counterclockwise order so that the connected component
$K_{-1}\ni -1$ of $K(P)\setminus\{\zeta\}$ lies between $r_1$ and $r_2$;
the component $K_{+1}\ni +1$ then lies between $r_{m+1}$ and $r_{m+2}$.
We label the remaining components of $K(P)\setminus \{\zeta\}$ as
$K_{\pm 1}^1,\ldots,K_{\pm 1}^{q-1}$, so that 
$K_{\pm 1}^i\ni P^{\circ i}(\pm 1)$. The corresponding objects associated
to $Q$ are similarly denoted with an added tilde. 

\realfig{pull-back-fig}{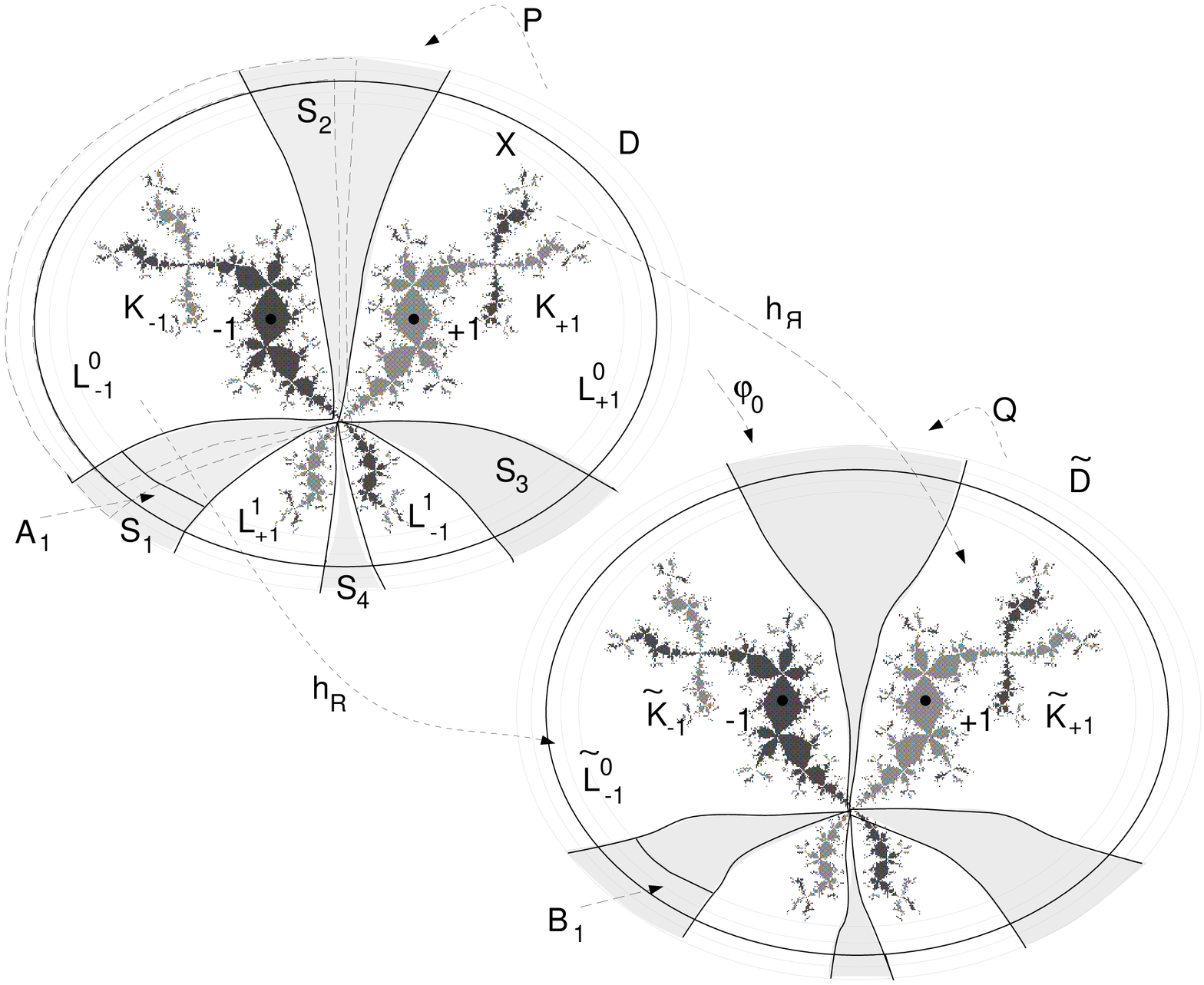}{Construction of the map $\varphi_0$ 
in the case $p/q=1/2$ with $m=3$}{11cm}

It will be convenient to introduce further notation.
Let $S_i\subset D$ be disjoint  invariant sectors centered at $r_i$,
and let $L^i_{\pm 1}$ be 
the component of $P^{\circ -q}(D)\setminus (\bigcup_{j=1}^{2q} S_j)$
containing $K_{\pm 1}^i$. 
The thickening procedure yields left and
right quadratic-like restrictions 
$$P^{\circ q}:\Omega'_{\lren}\to\Omega_{\lren}\;\;\n\;\;P^{\circ q}:
\Omega'_{\rren}\to\Omega_{\rren}$$
and
$$Q^{\circ q}:\tl\Omega'_{\lren}\to\tl\Omega_{\lren}\;\;\n\;\;Q^{\circ q}:
\tl\Omega'_{\rren}\to\tl\Omega_{\rren}.$$ 
By assumption, there exist   hybrid equivalences 
$h_{\lren}$ between $P^{\circ q}|_{\Omega'_{\lren}}$ and 
$Q^{\circ q}|_{\tl\Omega'_{\lren}}$, and  $h_{\rren}$ between 
$P^{\circ q}|_{\Omega'_{\rren}}$ and $Q^{\circ q}|_{\tl\Omega'_{\rren}}$.
We now replace the domains $D$ and $\tl D$ by $X=P^{\circ -q}(D)$ 
and $\tl X=Q^{\circ -q}(\tl D)$. 
We define the map $\varphi_0$ on $\bigcup_{i=0}^{q-1}L_{\pm 1}^i\subset X$ as
$$\varphi_0(z)=\left\{\begin{array}{lll}
		h_{\lren}(z) & \mbox{ for } & z\in L_{-1}^0,\\
		h_{\rren}(z) & \mbox{ for } & z\in L_{+1}^0,\\
		Q^{\circ -(q-i)}\circ h_{\lren}\circ P^{\circ q-i}(z) &
                    \mbox{ for } & z\in L_{-1}^i \mbox{ with } i>0,\\
		Q^{\circ -(q-i)}\circ h_{\rren} \circ P^{\circ q-i}(z) & 
	            \mbox{ for } & z\in L_{+1}^i \mbox{ with } i>0.\\
		\end{array}\right.$$
where by $Q^{\circ -(q-i)}$ we understand the univalent branch mapping
$\tl K_{\pm 1}$ to $\tl K_{\pm 1}^i$. 


Let $A_i$ be the strip of $X\setminus P^{\circ -q}(X)$ contained in $S_i$, and
$B_i$ its counterpart in $\tl X\setminus Q^{\circ -q}(\tl X)$.
We smoothly extend $\varphi_0$ to map $A_i\to B_i$ in agreement
with the previously specified values of $\varphi_0$ on
$\partial(L_{-1}^i\cup L_{+1}^i)$ and 
$\varphi_0\circ P^{\circ q}=Q^{\circ q}\circ \varphi_0$ on the inner
boundary of $A_i$. We now extend $\varphi_0$ to the entire sector
$S_i\cap X$ by setting 
$\varphi_0(z)=Q^{\circ -nq}\circ\varphi_0\circ P^{\circ nq}(z)$ when
$P^{\circ nq}(z)\in A_i$. 

The quasiconformal homeomorphism $\varphi_0:X\to\tl X$ so defined conjugates 
$P$ to $Q$ on the set 
$\bigcup_{i=0}^{q-1} P^{\circ i}(K_{\lren}\cup K_{\rren})$, 
with $\bar\partial \varphi_0=0$ almost everywhere on this set, sending each
$P^{\circ i}(K_{\lren})$ to $Q^{\circ i}(\tl K_{\lren})$ and each
$P^{\circ i}(K_{\rren})$ to $Q^{\circ i}(\tl K_{\rren})$.
We further extend $\varphi_0$ to a quasiconformal homeomorphism from
$Y=P(X)$ to $\tl Y=Q(\tl X)$ so that 
$$\varphi_0\circ P|_{\partial X}=Q\circ \varphi_0|_{\partial X}.$$
As $\varphi_0$ is a conjugacy between postcritical sets, there
is a unique lift $\varphi_1:X\to \tl X$ agreeing with $\varphi_0$ on
$\bigcup_{i=0}^{q-1} P^{\circ i}(K_{\lren}\cup K_{\rren})$ such that
the following diagram commutes:
$$\begin{CD}
X @>\varphi_1>> \tl X\\
@V{P}VV  @V{Q}VV\\
Y @>\varphi_0>> \tl Y
\end{CD}$$

As in the proof of \propref{inverse}, we set $\varphi_1(z)=\varphi_0(z)$
for $z$ in the
annulus $Y\setminus X$, and iterate the lifting procedure to obtain 
a sequence of quasiconformal maps $\varphi_n$ with uniformly bounded
dilatation. In view of the density of $\check{K}(P)$ in $K(P)$, the limiting
map $\varphi:Y\to \tl Y$ conjugates $P$ to $Q$.
As $\varphi_n$ stabilizes pointwise on $\check{K}(P)$ with
$\varphi|_{K_{\lren}}=h_{\lren}$ and $\varphi|_{K_{\rren}}=h_{\rren}$
by construction, it follows from Bers' \lemref{rickman} that
$\bar\partial \varphi=0$ almost everywhere on $\check{K}(P)$.
\end{pf}

\medskip
To conclude the proof of \propref{renorm-inj}, we show that the conjugacy 
just obtained is actually a hybrid equivalence: that any 
measurable invariant linefield on $K(P)\setminus\check K(P)$ has support in
a set of Lebesgue measure 0. In view of \lemref{pull-back}, it follows from  
the standard considerations of 
parameter dependence in the Measurable Riemann Mapping Theorem (\ref{MRMT})
that ${\cal F}=(\lrren)^{-1}(A,D)$
is the injective complex-analytic image of a polydisc ${\BBB D}^k$ for some
$k\in\{0,1,2\}$; see \cite{MSS} and \cite{McS}. On the other hand,
$\cal F$ is compact by Proposition \ref{proper}, whence $k=0$ and $\cal F$ 
is a single point.  $\;\;\Box$

\subsection{Conclusions}
\label{conclusions}

Setting $$h_{p/q,m}(c,\tl c)=(\lrren)^{-1}(\kappa_{p/q}(c),
\kappa_{p/q}(\tl c))$$
so that
$h_{p/q,m}(c,\tl c)=(A,D)$ if and only if  $P_{A,D}=f_c
{\underset{p/q,m}{\intr}} f_{\tl c}$
we obtain the embeddings
$$h_{p/q,m}:(\M_{p/q}\setminus\{\rt_{p/q}\})\times(\M_{p/q}\setminus
\{\rt_{p/q}\})\longrightarrow B_{p/q,m}\subset C_{p/q,m}$$
whose existence was asserted in \thmref{main}.
As observed in \secref{outline}, if $P=P_{A,D}$ is birenormalizable
and $(A,D)\in\CubMR$ then $p/q=0$ or $1/2$. It follows by symmetry that 
$\lrren(A,D)=(c,\bar c)$ for some $c\in\BBB C$; conversely, if
$\lrren(A,D)=(c,\bar c)$ then $(A,D)\in\CubMR$.
Writing $\bar\Delta$ for the antidiagonal
embedding $\BBB C\ni c\mapsto(c,\bar c)\in{\BBB C}^2$, we define
$$\Psi_{1/2,1}= h_{1/2,1}\circ \bar\Delta
:\M_{1/2}\setminus\{\rt_{1/2}\}\to\CubMR\cap C_{1/2,1}$$
$$\Psi_{1/2,3}= h_{1/2,3}\circ \bar\Delta
:\M_{1/2}\setminus\{\rt_{1/2}\}\to\CubMR\cap C_{1/2,3}$$
and
$$\Psi_0= h_0\circ \bar\Delta:\M\setminus\{\rt\}\to\CubMR\cap C_0.$$
These are the embeddings whose existence was asserted
in \thmref{real-main-thm}. Compatibility with the standard
planar embeddings is a consequence of the following recent result of
Buff \cite{Buff}:

\begin{thm}
\label{ext}
Let $K_1$ and $K_2$ be compact, connected, cellular sets in the plane, and
$\varphi:K_1\to K_2$ a homeomorphism. If $\varphi$ admits a continuous
extension to an open neighborhood of $K_1$ such that points outside $K_1$
map to points outside $K_2$, then $\varphi$ extends to a homeomorphism
between open neighborhoods of $K_1$ and $K_2$. 
\end{thm}

Let us sketch the argument for the map $\Psi_0$. It is easily verified from
the explicit expressions in \cite[p. 22]{Mil} that for each $\mu\in{\BBB C}
\setminus\{1\}$ there is a unique pair $(A,B)\in\BBB R^2$ such that the
corresponding polynomial in the normal form (\ref{normal-form}) has a pair of
complex conjugate fixed points with multipliers $\mu$ and $\bar\mu$, the
remaining fixed point having eigenvalue
$$\nu=1+\frac{|\mu-1|^2}{2\Re(\mu-1)}.$$ 
We may continuously label these multipliers as $\mu(A,D)$, $\bar\mu(A,D)$
and $\nu(A,D)$ for parameter values $(A,D)\in\BBB R\times i\BBB R$
in a neighborhood of $\Psi_0(\M\setminus\{\rt\})$;
in particular, $(A,D)\mapsto\mu(A,D)$ is a homeomorphism on such a 
neighborhood. It follows from Yoccoz Inequality (\ref{yoccoz}) that
$\nu(A,D)>1$, and therefore $\mu(A,D)\not\in[1,\infty)$, for $(A,D)$ in 
$\Psi_0(\M\setminus\{\rt\})$. Similarly, $\mu(\Psi_0(c))\to 1$ as 
$c\to\rt$, and thus $c\mapsto\mu(\Psi_0(c))$ extends to a embedding
$$\Upsilon:\M\to \BBB C\setminus(1,\infty)$$ which clearly commutes
with complex conjugation.

We claim that $\Upsilon^{-1}:\Upsilon(\M)\to \M$
admits a continuous extension meeting the condition of \thmref{ext}.
The idea is to allow renormalizations with disconnected Julia sets.
Recalling \lemref{same-rays}, we note that
the rays $r_0$ and $r_{1/2}$ continue to land at the same fixed point
for $(A,D)$ in a neighborhood of $\Psi_0(\M\setminus\{\rt\})$. 
As before, we may construct left and right quadratic-like restrictions 
with continuously varying domains $\Omega'_{A,D}$.
It is emphasized in Douady and Hubbard's
original presentation \cite{DH} that straightening, while no longer canonical
for maps with disconnected Julia set, may still be continuously defined: 
it is only necessary to begin with continuously varying 
quasiconformal homeomorphisms from the fundamental annuli 
$\Omega_{A,D}\setminus\Omega'_{A,D}$ to the standard annulus. We
thereby obtain a continuous extension to a neighborhood of
$\Upsilon(\M\setminus\{\rt\})$; it is easily arranged that this extension
commutes with complex conjugation, so that it is trivial to obtain a
further extension to an open set containing the point $1$.

\section{Measure of the Residual Julia Set}
\label{res-mes}

Recall that for a birenormalizable polynomial $P$,
$$\check{K}(P)=\bigcup_{i=0}^{\infty}P^{-\circ i}(K_\lren\cup K_{\rren}).$$
Here we synthesize various arguments of Lyubich to show that the
{\it residual Julia set} $K(P)\setminus\check{K}({P})$ 
has Lebesgue measure 0, provided that neither 
renormalization lies in the closure of the main hyperbolic component of
$\M$. Subject to this restriction, we arrive at an alternative proof that 
the conjugacy constructed in \lemref{pull-back} is a hybrid
equivalence. We formalize the statement as follows:

\begin{thm}
\label{rel-meas-0}
Let $P$ be a birenormalizable cubic polynomial. If all three fixed
points are repelling then $K(P)\setminus\check{K}(P)$ has Lebesgue measure 0.
\end{thm}

The main technical tool for us will be the celebrated Yoccoz puzzle 
construction which we briefly recall below:

\medskip

\noindent {\bf Yoccoz puzzle and recurrence.}
Let $f=f_c$ be a quadratic polynomial with connected Julia set, and 
$G\supset K(f)$ be a domain bounded by some fixed equipotential curve.
As observed in \lemref{rot-number}, if $c\in L_{p/q}$ for some $q\geq 2$ then
$\alpha_f$ is the landing point of a cycle of $q$ external
rays. The {\it Yoccoz puzzle of depth zero} consists of the $q$ pieces 
$Y^0_1,$ $Y^0_2,\ldots, Y^0_q$ obtained by cutting 
$G$ along these rays, and  
the puzzle pieces of depth $n>0$ are the connected components $Y^n_i$
of the various $f^{\circ -n}(Y^0_i)$.
Each point $z\in K(f)\setminus f^{\circ -n}(\alpha)$ 
lies in a unique depth $n$ puzzle piece $Y^n(z)$.
A nonrenormalizable 
 polynomial $f$ has a {\it reluctantly recurrent} critical point
if there exists $k\geq 0$ and  a sequence of depths $n_i\to \infty$ such that
the restriction $f^{\circ n_i-k}:Y^{n_i}(0)\to Y^k(f^{\circ n_i-k}(0))$
has degree $2$. Note that, somewhat abusing the notation, we allow maps
with non-recurrent critical point in this definition. In the complementary
case of {\it persistently recurrent} critical point 
Lyubich has shown the following:

\comm{
For fixed nonrenormalizable $f=f_c$, the 
Yoccoz $\tau$-function is defined as follows. Given $n\in \BBB N$
we set $\tau(n)$ equal to the largest $m\leq n$ such that the
depth $m$ puzzle-piece 
$Y^m(f^{\circ n-m}(c))=f^{\circ n-m}(Y^n(c))$
contains the critical point $0$;
by convention, $\tau(n)=-1$ if no such $m$ exists.
We say that $f$ has a {\it persistently recurrent} critical point
if $\lim_{n\to\infty}\tau(n)=\infty$. 
}

\begin{lem}{\cite[p. 6]{Ly1}}
\label{minimal}
If the critical point of a non-renormalizable quadratic polynomial
$f_c$ is persistently recurrent then $f_c|_{\omega_c(0)}$ is
topologically minimal, that is all orbits are dense in $\omega_c(0)$.
In particular, $\beta_{f_c}\not\in\omega_c(0)$.
\end{lem}
   
The puzzle construction is easily adapted to a cubic map  $P=P_{A,D}$ 
which has a connected Julia set with empty interior 
and every periodic orbit repelling.
The depth zero puzzle pieces $W^0_i$ are now obtained by
cutting an equipotentially bounded domain $G\supset J(P)$ 
along every ray which lands at some fixed point, and
the pieces of depth $n$ are the connected components
of the various $P^{\circ -n}(W^0_i)$. Each point $z\in J(P)\setminus
P^{-\circ n}(\{\text{fixed points}\})$ lies in a unique  depth $n$ puzzle
piece $W^n(z)$.

By analogy with the quadratic case, we say that the critical point
$\pm 1$ of the cubic polynomial $P=P_{A,D}$ is {\it reluctantly recurrent} 
if there exist $k\geq 0$, $N>0$
 and a sequence of depths $n_i\rightarrow\infty$
such that $P^{\circ n_i-k}|_{W^{n_i}(\pm 1)}$ is a map of degree $N$.
We readily observe that if $P$ is birenormalizable and one of its
renormalizations has a reluctantly recurrent critical point
then the corresponding critical point of $P$ is reluctantly recurrent.
Indeed in this case the restriction 
$P^{\circ n-k}|_{W^{n}( 1)}$ has the same degree as the map 
$\lren(P)^{\circ n-k}$ on the quadratic puzzle piece $Y^n(0)$,
 and similarly for the other renormalization.

Yarrington \cite{Yar} has shown that if both critical points of $P$ are
reluctantly recurrent then $J(P)$ is locally connected;  in particular
nested sequences of puzzle pieces shrink to points in this case:

\begin{equation}\label{shrink}
\bigcap_{n=0}^{\infty} W^n(z)=\{z\}
\end{equation}
for every $z\in J(P)\setminus \{\text{prefixed points}\}$ 
(\cite[Theorem 3.5.7]{Yar}).

\medskip
\noindent{\bf Relative ergodicity. }
The proof of \thmref{rel-meas-0} is based on the following general principle
of Lyubich:

\begin{thm}[\cite{Ly}]\label{asympt}
Let $g$ be a rational map with $J(g)\neq\hat{\BBB C}$. Then
$$\omega_g(z)\subset \bigcup_{\gamma\in\Gamma}\omega_g(\gamma)$$
for almost every $z\in J(g)$, where $\Gamma$ is the set of all critical points.
\end{thm}

We divide the argument into two cases depending on the recurrence properties
of the renormalizations of $P$.

Assume first that $\omega_P(-1)\cap \omega_P(+1)=\emptyset$.
It follows from \thmref{asympt} that for almost every $z\in J(P)$
$\omega_P(z)$ lies in the disjoint union $\omega_P(-1)\cup\omega_P(+1)$.  
Without loss of generality $\omega_P(z)\subset\omega_P(-1)$, so
there is a subsequence $P^{\circ nq+\ell}(z)$ such that every accumulation
point lies in $K_{\lren}$. In particular, $P^{\circ nq+\ell}(z)\in\Omega'$ 
for sufficiently large $n$, where $\Omega'$ is the domain of the 
left quadratic-like restriction of $P^{\circ q}$. Thus, $P^{\circ nq+\ell}(z)
\in K_{\lren}$ for large enough $n$, and therefore
$z\in \bigcup_{i=0}^{\infty}P^{\circ -i}(K_{\lren})\subset\check{K}(P)$.

In the other case, recall from (\ref{bi-omega}) that 
$\omega_P(-1)\cap \omega_P(+1)=\{\zeta\}$.
Combining Remark \ref{not-accumulate} and 
\lemref{minimal}, we see that both  $\lren(P)$ and $\rren(P)$ 
are nonrenormalizable quadratics with reluctantly recurrent
critical points.  We conclude the argument  by showing that under
 these conditions the Lebesgue
measure of the Julia set of $P$ is zero:

\begin{lem}
\label{meas-0}
Let $P=P_{A,D}$ be a birenormalizable cubic with all
three fixed points repelling.
If both its renormalizations $\lren(P)$ and $\rren(P)$ are
nonrenormalizable quadratic maps  with reluctantly recurrent 
critical points, then the Julia set of $P$ has Lebesgue measure zero.
\end{lem}
\begin{pf}
We  adapt Lyubich's argument \cite{Ly1} for the quadratic case.
As both  critical points of $P$ are reluctantly recurrent,
there exist $k$ and arbitrary large $s$ and $t$ such that
$P^{\circ s}|_{W^{s+k}(-1)}$ and $P^{\circ t}|_{W^{t+k}(+1)}$ are maps of
degree 2. By \thmref{asympt} for a full measure set of $z\in J(P)$ 
there exists $n$ such that $P^{\circ n}(z)$ lies in
$W^{s+k}(-1)\cup W^{t+k}(+1)$ for any $s$ and $t$.
Fixing $t$, $s$ and  $z$ consider the least such $n$. 
Without loss of generality,
$P^{\circ n}(z)\in W^{s+k}(-1)$ and we obtain a chain of univalent branches
of $P^{-1}$ 
$$W^{s+k}(-1)=X_0\leftarrow X_{-1}\leftarrow\ldots\leftarrow X_{-n}\ni z$$ 
by pulling this piece back along the orbit of $z$.

Fix a puzzle-piece $W_i^k$  of depth $k$. As the boundary of 
$W^k_i$ consists of preimages of external rays landing at the 
fixed points of $P$ and equipotential curves it follows from Koebe
$1/4$ theorem that there exist
$\delta_i$ and $\delta_i'$ such that for any $u\in W^k_i$
with $\dist(\partial W^k_i,u)<\delta_i$  some neighborhood 
$U\subset W^k_i$ around $u$  is univalently mapped by an iterate
$P^{\circ j}$ to a disk  $D_{\delta'_i}(P^{\circ j}(u))$ of radius $\delta_i'$ 
centered at $P^{\circ j}(u)$.
Denote by $\delta$ the minumum of various $\delta_i$, $\delta_i'$
and set $u=P^{\circ s+n}(z)
\in P^{\circ s}({W^{s+k}(-1)})=W^k$. By the above, there 
exists a neighborhood $U\ni u$ in $W^k$ and an iterate  $P^{\circ j}$ 
univalently mapping $U$ to $D_\delta(P^{\circ j}(u))$. 

\smallskip
Assume first that $P^{\circ s}(-1)$ does not belong to $U$, or
 $|P^{\circ s+j}(-1)-P^{\circ j}(u)|>\delta/100$. The density of the 
Julia set in  a disk of radius $\delta/200$ is bounded away from $1$.
Consider the univalent pullback $T_0=D_{\delta/200}(P^{\circ j}(u)),
T_{-1},\ldots,T_{-s-n-j}$ along the orbit $z\mapsto P(z)\mapsto\ldots
\mapsto P^{\circ s+n+j}(z)\in T_0$. By the Koebe distortion theorem,
the density of $J(P)$ in $T_{-s-n-j}$ is also bounded away from  $1$. 
By the estimate (\ref{shrink}), the disks $T_{-s-n-j}$ shrink to
the point $z$ as $s$ grows  and therefore $z$ is not a point of density
for the set $J(P)$. 

\smallskip
Consider now the case when
$P^{\circ s}(-1)\in U$, and $|P^{\circ s+j}(-1)-P^{\circ j}(u)|<\delta/100$.
Then we can find a disk $D_1$ centered at $P(-1)$, such that
$P^{\circ s+j-1}(D_1)$ is contained between $D_{\delta/10}(P^{\circ j}(u))$
and $D_{\delta/2}(P^{\circ j}(u))$. By Koebe distortion theorem,
the density of $J(P)$ in the disk $D_1$ is bounded away from $1$. 
Consider the preimage $D_0$ of $D_1$ centered around $-1$ and contained
in $W^{n+s}(-1)$. The density of the Julia set in $D_0$ is again bounded
away from $1$, and  as in the previous case we conclude that $z$ is not a 
point of density.
Thus the set of density points of $J(P)$ has measure zero, and by Lebesgue
density theorem so does $J(P)$.
\end{pf}

\appendix

\section{Discontinuity at the corner point}
\label{discont-root}

The systematic exclusion of root points is not merely an artifact of
our reliance on quasiconformal surgery. It is conceivable that more powerful
techniques might someday prove existence and uniqueness of 
intertwinings $f_c{\underset{p/q,m}{\intr}}f_{\tl c}$
for {\it any} pair \linebreak $(c,\tl c)\in\M_{p/q}\times\M_{p/q}$.
Indeed, $f_c|_{J(f_c)}$ for any $c\in H_{p/q}$ is canonically topologically
conjugate to $f_{\hat{c}}|_{J(f_{\hat{c}})}$ where $\hat{c}=\rt_{p/q}$, and on 
these grounds we have put forth Conjecture \ref{wishful-thinking}.
On the other hand, we assert in \thmref{two-disks} that
such a extension of $h_{p/q,m}$ is necessarily discontinuous
at the corner point $(\rt_{p/q},\rt_{p/q})$. This is an
instance of a phenomenon investigated by one of the authors. It is 
shown in \cite{mate} that any disjoint type component consisting of
maps with {\it adjacent} attracting basins must suffer such a discontinuity;
for $c,\tl c\in H_{p/q}$ the basins of $f_c{\underset{p/q.m}{\intr}} 
f_{\tl c}$ are adjacent by construction.
Here we simply summarize the relevant considerations.

Let $g$ be an analytic map fixing $\zeta\in \BBB C$. The
{\it holomorphic index} of $g$ at $\zeta$ is the residue
$$\eta=\frac{1}{2\pi i}\int_\gamma \frac{dz}{z-g(z)}$$
where $\gamma$ is a loop enclosing $\zeta$ but no other
other fixed point. It is easily checked that this quantity is conformally
invariant; in fact, 
$$\eta=\frac{1}{1-\lambda}$$
so long as the multiplier $\lambda=g'(\zeta)$ is not equal to 1. 

An elementary computation yields
$$\eta_a=\frac{1}{2\pi i}\int_\gamma\frac{dz}{-az^2-z^3}=\frac{1}{a^2}$$
for the holomorphic index of $Q_a(z)=z+az^2+z^3$ at the parabolic fixed
point $0$. In the terminology of \cite{mate}, such a fixed point is
described as {\it parabolic-attracting}, {\it parabolic-indifferent}
or {\it parabolic-repelling} depending on whether $\Re \eta$ is greater than,
equal to, or less than 1. The first of these alternatives applies when
$|a^2-\frac{1}{2}|<\frac{1}{2}$. The corresponding region in the
$a-$plane is bounded by a lemniscate shaped like the symbol $\infty$; its
position in the cubic connectedness locus is depicted in \figref{per11}, which
in view of the 4-dimensionality of $D_{1,1}$ is merely schematic. The
intersection of the component boundary with $\Per_1(1)$ consists of the
closure of this lemniscate (shaded in dark gray, and contained in the light
gray region where both critical points lie in the parabolic basin), and a 
similar locus (the large lobes of the medium gray region) parameterizing maps
whose other fixed point is attracting or indifferent; the latter might
be described as intertwinings $f_c{\underset{0}{\intr}} 
f_{\tl c}$ for $c\in H_0$ and $\tl c=\rt$. These pieces 
intersect at the parameter value $a=0$ where the parabolic fixed point 
becomes {\em degenerate}.
	
\realfig{per11}{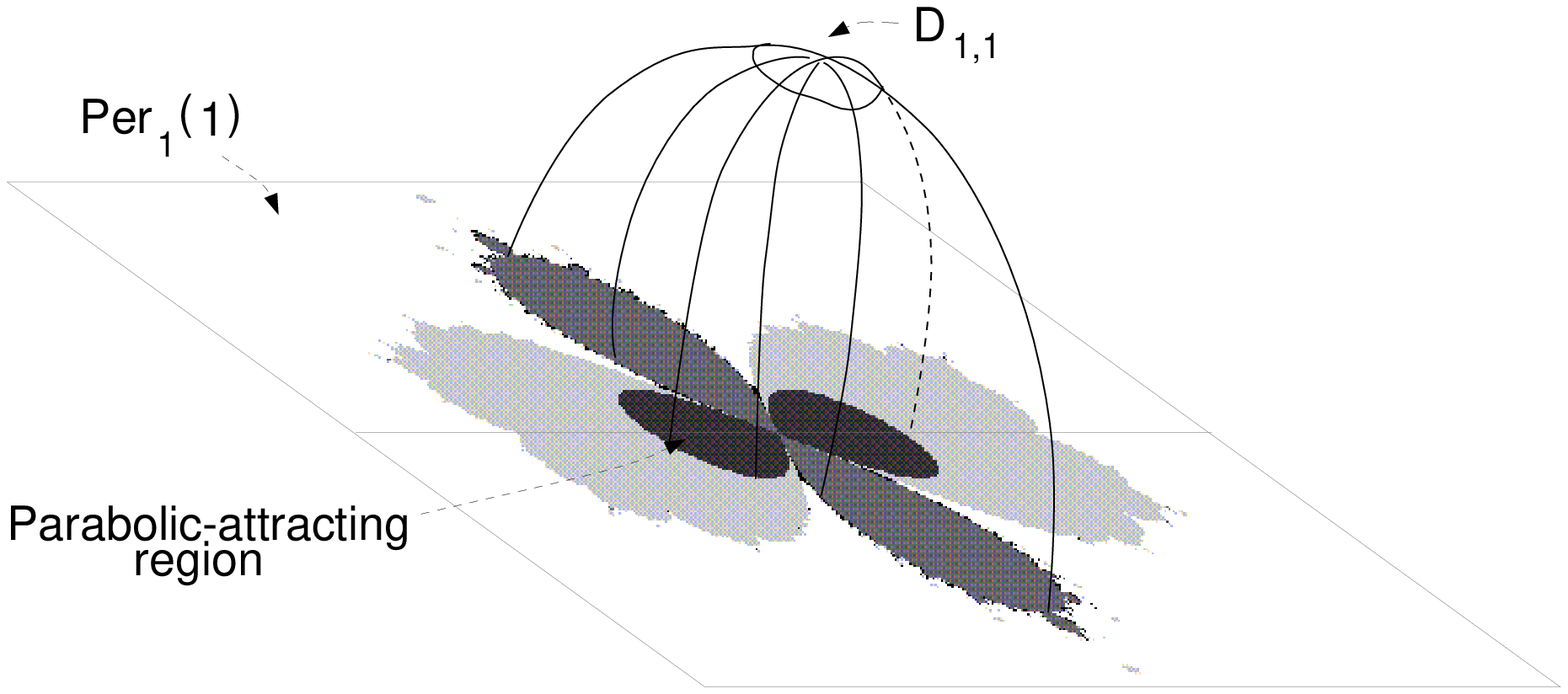}{The unique component of type $D_{1,1}$ and its
impression on $\Per_1(1)$}{10cm}

The crux of the matter is the following elementary observation 
(compare with \cite[Problem 9-1]{Mil-survey}):
\begin{lem}
\label{perturb}
Let $\eta\in\BBB C$. Then $\,\Re\eta\geq 1$ if and only if there exist
continuous paths $\lambda,\tl\lambda:[0,1)\to\BBB D$ with endpoints
$\lambda(1)=1=\tl\lambda(1)$ such that 
$$\eta=\lim_{t\to 1}\;\frac{1}{1-\lambda(t)}+\frac{1}{1-\tl\lambda(t)}.$$
Complex conjugate paths may be chosen when $\eta$ is real.
\end{lem}

Similar considerations apply to $h_{p/q,m}(H_{p/q}\times H_{p/q})$. For odd 
denominator $p/q$ and $A_{p/q}=-\frac{1}{3}e^{2\pi ip/q}$, it is
easy to see that $P_{A_{p/q},0}$ is the unique normalized cubic polynomial
with a degenerate parabolic fixed point of multiplier $e^{2\pi ip/q}$; thus
$$\lim_{c\rightarrow\rt_{p/q}} h_{p/q,q}(c,c)\to (A_{p/q},0)$$
as is evident in \figref{symmlocus}.

\section{Non local-connectivity of the real connectedness locus}
\label{section-comb}

Here we employ a simplified version of an argument of Lavaurs \cite{Lav}
to conclude that the real cubic connectedness locus is not locally
connected along an interval in the boundary of $\Psi_0(\M\setminus\{\rt\})$.
The existence of comb-like structures in \linebreak
$\Psi_{1/2,i}(\M_{1/2}\setminus\{\rt_{1/2}\})$ is similarly demonstrated; 
see also
Nakane and  Schleicher's proof of non local-connectivity for the
tricorn \cite{Nakane}.

We begin with a brief review of the theory of parabolic bifurcations, as
applied in particular to real cubic polynomials.
The reader is referred to \cite{Do} for a more comprehensive exposition;
supporting technical details may be found in \cite{Sh}. Recall that the
fixed point at 0 is parabolic with multiplier 1 for every map in the family
$$Q_a(z)=z+az^2+z^3.$$

\begin{lem}[Fatou coordinates]
\label{petals}
For $a\neq 0$ there exist topological discs $U^\ap_a$ and $U^\rp_a$
whose union is a punctured neighborhood of the parabolic fixed point, such that
\begin{eqnarray*}
Q_a(\bar U^\ap_a)\subset U^\ap_a\cup\{0\} & \mbox{ and } & 
\bigcap_{k=0}^\infty Q_a^{\circ k}(\bar U^\ap_a)=\{0\},\\
Q_a^{-1}(\bar U^\rp_a)\subset U^\rp_a\cup \{0\} & \mbox{ and } & 
\bigcap_{k=0}^\infty Q_a^{\circ -k}(\bar U^\rp_a)=\{0\}.
\end{eqnarray*}
Moreover, there exist injective analytic maps
$$\Phi^\ap_a:U^\ap_a\to \BBB C\;\; \n\;\; \Phi^\rp_a:U^\rp_a\to \BBB C,$$
unique up to post-composition by translations, such that
$$\Phi^\ap_a(Q_a(z))=\Phi^\ap_a(z)+1\;\; \n\;\;
\Phi^\rp_a(Q_a(z))=\Phi^\rp_a(z)+1.$$
The quotients $C^\ap_a=U^\ap_a/Q_a$ and 
$C^\rp_a= U^\rp_a/Q_a$ are therefore Riemann surfaces
conformally equivalent to the cylinder $\BBB C/\BBB Z$.
\end{lem}

\medskip
The quotients $C^\ap_a$ and $C^\rp_a$ are customarily referred to as the 
{\it \'Ecalle-Voronin cylinders} associated to the map $Q_a$; we will find
useful to regard these as Riemann spheres with distinguished
points $\pm$ filling in the punctures. Every point in the {\it parabolic basin}
$${\cal B}_a=\{z\in\BBB C|\; Q_a^{\circ n}(z)\neq 0 \mbox{ for } n\geq 0 
\mbox{ and } Q_a^{\circ n}(z)\to 0\}$$
eventually lands in $U^\ap_a$, and the return map from $U^\rp_a\cap{\cal B}_a$
to $U^\ap_a$ descends to a well-defined analytic transformation 
$${\cal E}_a:{\cal W}_a\to C^\ap_a$$
where ${\cal W}_a$ is the image of ${\cal B}_a$ on $C^\rp_a$. It is easy to
see that the ends of $C^\rp_a$ belong to different components of 
${\cal W}_a$. The choice of a conformal
{\em transit isomorphism} $$\Theta:C^\ap_a\to C^\rp_a$$ respecting these ends
determines an analytic dynamical system
$${\cal F}_{a:\Theta}=\Theta\circ{\cal E}_a:{\cal W}_a\to C^\rp_a$$
with fixed points at $\pm$.  The product of the corresponding
eigenvalues $\varrho^\pm_{a:\Theta}$ is clearly independent of $\Theta$,
indeed $$\varrho^+_{a:\Theta}\cdot\varrho^-_{a:\Theta}=e^{-4\pi^2(\eta_a-1)}
=e^{-4\pi^2(\frac{1}{a^2}-1)}.$$

For $a\in\BBB R$ the real-axis projects to natural {\it equators}
$R^\ap_a\subset C^\ap_a$ and $R^\rp_a\subset C^\rp_a$; the set
${\cal W}_a$ is disjoint from $R^\rp_a$ and symmetric about it.
Moreover, when $a\in(0,\sqrt{3})$ the critical points of $Q_a$ form a
complex conjugate pair in ${\cal B}_a$. We restrict attention to this
simplest case: $J(Q_a)=\partial{\cal B}_a$ is a Jordan curve, as is each of
the two components of $\partial{\cal W}_a$. It follows from the details of
the construction that ${\cal E}_a$ has infinitely many critical points but
only two critical values; these are situated symmetrically with respect
to the appropriate equators, and each of the critical values $\upsilon^\pm_a$
has critical preimages on both sides of $R^\rp_a$. 

We now consider perturbations in the family 
$$Q_{a,\eps}(z)=\eps+z+az^2+z^3$$
corresponding to
\begin{equation} \label{BHQ}
F(z)=z^3-\frac{1}{3}a^2z+\left(\frac{2}{27}a^3+\eps\right)
\end{equation}
in the normal form of (\ref{normal-form}).
For small $\eps>0$ the parabolic point splits into a complex conjugate
pair of attracting fixed points $\zeta^\pm_{a,\eps}$ but one may still
speak of attracting and repelling petals:

\begin{lem}[Douady coordinates]
\label{pert-coord}
For small $\eps>0$ there exist topological discs $U^\ap_{a,\eps}$ and
$U^\rp_{a,\eps}$ whose union is a neighborhood of the parabolic fixed point
of $Q_a$, and injective analytic maps
$$\Phi^\ap_{a,\eps}:U^\ap_{a,\eps}\to \BBB C\;\; \n\;\;
\Phi^\rp_{a,\eps}:U^\rp_{a,\eps}\to \BBB C,$$
unique up to post-composition by translations, such that
$$\Phi^\ap_{a,\eps}(Q_{a,\eps}(z))=\Phi^\ap_{a,\eps}(z)+1\;\;
\n\;\; \Phi^\rp_{a,\eps}(Q_{a,\eps}(z))=\Phi^\rp_{a,\eps}(z)+1.$$
The quotients $C^\ap_{a,\eps}= U^\ap_{a,\eps}/Q_{a,\eps}$
and $C^\rp_{a,\eps}= U^\rp_{a,\eps}/Q_{a,\eps}$ are Riemann surfaces
conformally equivalent to $\BBB C/\BBB Z$.
\end{lem}

In view of the assumption on $\eps$ these cylinders come similarly equipped
with equators.  As in the parabolic case, the return map
from the relevant portion of $U^\rp_{a,\eps}$ to 
$U^\ap_{a,\eps}$ descends to an analytic transformation ${\cal E}_{a,\eps}$ 
from a neighborhood of each end of $C^\rp_{a,\eps}$
to a neighborhood of the corresponding end of $C^\ap_{a,\eps}$. However, there
is now a canonical transit isomorphism $\Theta_{a,\eps}:C^\ap_{a,\eps}
\to C^\rp_{a,\eps}$,
and the composition $${\cal F}_{a,\eps}=\Theta_{a,\eps}\circ{\cal E}_{a,\eps}$$
is completely specified by the dynamics of $Q_{a,\eps}$. In particular, the
eigenvalues at $\pm$ are given by
\begin{equation}\label{eigprod}
\varrho^\pm_{a,\eps}=e^{\frac{-4\pi^2}{\log \lambda^\pm(a,\eps)}}
\end{equation}
where $\lambda^\pm(a,\eps)=1\pm 2i\sqrt{a\eps}+O(\eps)$ are the complex
conjugate eigenvalues of $\zeta^\pm_{a,\eps}$. 
A fixed but otherwise arbitrary choice of basepoints in the original petals
$U^\ap_{a,\eps}$ and $U^\rp_{a,\eps}$ allows us to identify $C^\ap_{a,\eps}$
with the various $C^\ap_{a,\eps}$ and $C^\rp_{a,\eps}$ with the various
$C^\rp_{a,\eps}$. The following 
fundamental theorem first appeared in \cite{orsay-notes} and was adapted
to the case at hand in \cite{Lav}:

\realfig{per11comb}{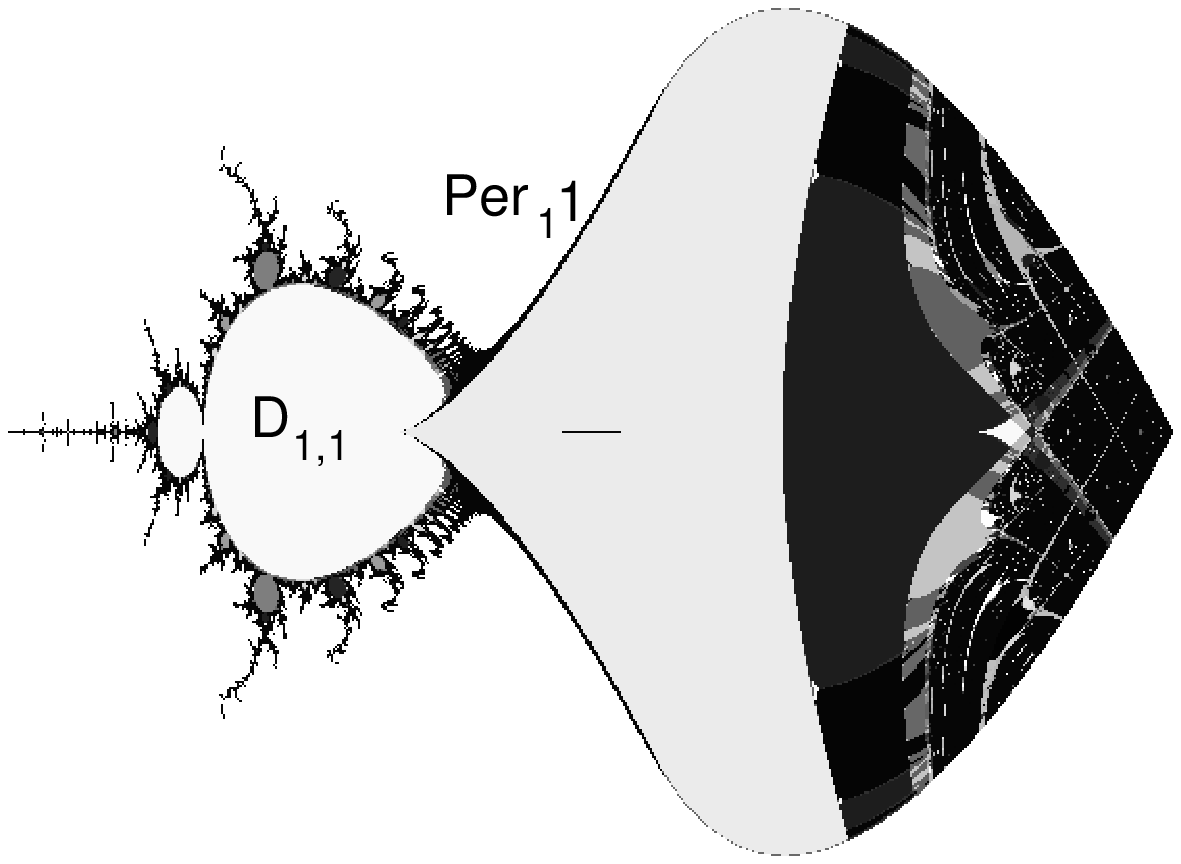}{The  view of the comb on the $D_{1,1}$ 
component in $(A,b)$ parametrization}{0.56\hsize}

\realfig{contfig}{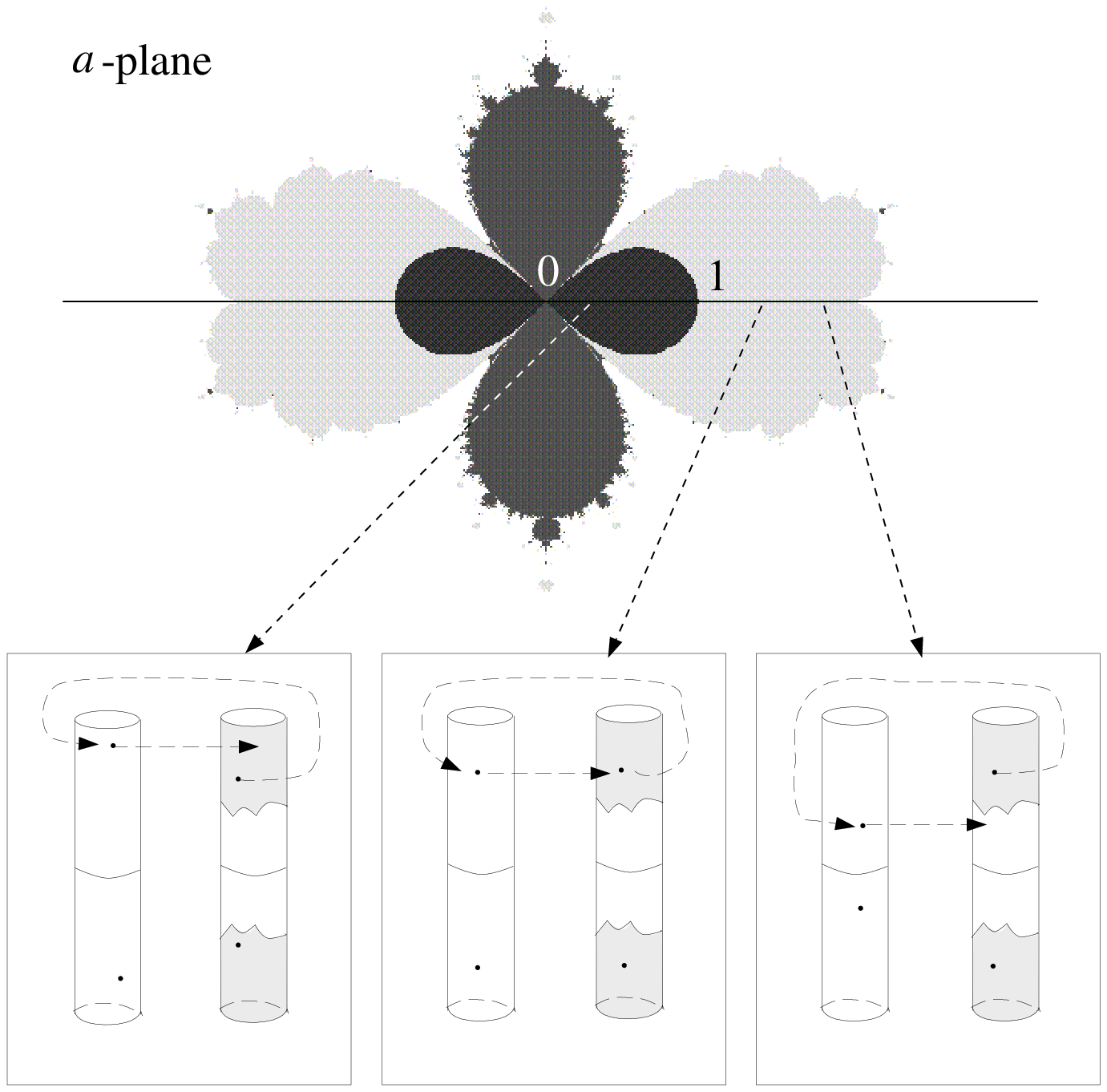}{Lavaurs' continuity argument}{0.6\hsize}

\begin{thm}
\label{conv-cyl}
In this setting, if $a_k\to a$ and $\eps_k\to 0$ such that 
$$\varrho^+_{a_k,\eps_k}\to\varrho^+_{a:\Theta}\;\mbox{ or equivalently }\;
\varrho^-_{a_k,\eps_k}\to\varrho^-_{a:\Theta}$$
then ${\cal F}_{a_k,\eps_k}\to{\cal F}_{a:\Theta}$
locally uniformly on ${\cal W}_a$.
\end{thm}

\medskip
We are now in a position to avail ourselves of an elementary but crucial
observation of Lavaurs \cite{Lav}:

\begin{lem} \label{conarg}
There exist $a\in(0,\sqrt{3})$ and a transit map $\Theta:C^\ap_a\to C^\rp_a$
respecting equators, such that both $\Theta(\upsilon^\pm_a)$ are 
superattracting fixed points for ${\cal F}_{a:\Theta}$.
\end{lem}

The relevant continuity argument is depicted in \figref{conarg}. For small
$a>0$, the critical values $\upsilon^+_a$ and $\upsilon^-_a$ are farther apart
than any pair of critical points of ${\cal E}_a$. All of these points move
continuously as $a$ increases towards the parameter value $\sqrt{3}$
where $\upsilon_a^\pm$ collide at the equator. Consequently, there
exist $a\in(0,\sqrt{3})$ and a symmetric pair of critical points $c_a^\pm$
which are exactly as far apart as the critical values $\upsilon_a^\pm$.
Both possibilities 
$${\cal E}_a(c_a^\pm)=\upsilon_a^\pm \;\;\n\;\; {\cal E}_a(c_a^\pm)=\upsilon
_a^\mp$$ may be so arranged. Choosing the former, we see that $\Theta(
\upsilon_a^\pm)=c_a^\pm$ for a suitable transit map respecting equators;
in particular, each of $c_a^\pm$ is a superattracting fixed point for
${\cal F}_{a:\Theta}$.

\medskip
Let $a$ be the parameter value so obtained.
In view of (\ref{eigprod}) there exist real $\eps_k$ decreasing to $0$ 
with  $\varrho^\pm_{a,\eps_k}\to\varrho^\pm_{a:\Theta}$.
It follows from \thmref{conv-cyl} that the nearby fixed points of
of ${\cal F}_{a,\eps_k}$ are attracting. Their lifts
generate a complex conjugate pair of attracting periodic
orbits in the original dynamical plane, and thus $J(Q_{a,\eps_k})$ 
is connected; moreover, $Q_{a,\eps_k}$ is birenormalizable 
as the critical orbits are separated by the real-axis. 
The two ways of marking the critical points of $Q_a$
yield parameters $(A_\infty,\pm D_\infty)\in\Per_1(1)$
and corresponding parameters $(A_k,\pm D_k)\in\Phi_0(\M-\{\rt\})$
associated to the perturbations $Q_{a,\eps_k}$. It follows from (\ref{BHQ}) 
that $A_\infty<0$, and thus $(A_\infty,\pm D_\infty)$ are the 
endpoints of an interval $I$ on the simple arc
$${\cal P}=\{(A,D)\in\Per_1(1)|\; A<0\}.$$ The entire impression
$${\cal I}=\{(A,D)\in\CubMR|\; \Psi_0(c_j)\to (A,D)\mbox{ for some } c_j\in
\M\setminus\{\rt\} \mbox{ with } c_j\to\rt\}$$
lies in $\cal P$ by Yoccoz Inequality (\ref{yoccoz}); thus $I\subset
{\cal I}$, as $\cal I$ is connected and
$$(A_\infty,\pm D_\infty)=\lim_{k\to\infty}(A_k,\pm D_k)\in{\cal I}$$
by construction. It follows from \lemref{same-rays} and the considerations of
\lemref{perturb} that $\CubMR$ is non-locally connected
at every $(A,D)\in{\cal I}$
for which $P_{A,D}$ has a parabolic-repelling fixed point.

\end{document}